\documentclass{gtpart}

\usepackage{geometry}            
\usepackage[parfill]{parskip}    
\usepackage{amsfonts}
\usepackage{epstopdf}
\usepackage[pdftex]{color,graphicx}
\usepackage{hyperref}
\usepackage[all]{xy}

\title{Wrinkled fibrations on near-symplectic manifolds}
\author{Yank\i\ Lekili}
\givenname{Yank\i}
\surname{Lekili}
\address{Department of Mathematics, M.I.T., Cambridge MA 02139, USA} 	
\email{lekili@math.mit.edu}


\volumenumber{}
\issuenumber{}
\publicationyear{}
\papernumber{}
\startpage{}
\endpage{}
\doi{}
\MR{}
\Zbl{}
\received{}
\revised{}
\accepted{}
\published{}
\publishedonline{}
\proposed{}
\seconded{}
\corresponding{}
\editor{}
\version{}

\newcommand{\f}[1]{\mathbb{#1}}
\newcommand{\del}{\partial}

\newtheorem{theorem}{Theorem}[section]

\newtheorem{definition}[theorem]{Definition}
\newtheorem{corollary}[theorem]{Corollary}

\theoremstyle{definition}

\begin{document}
\begin{abstract}
Motivated by the programmes initiated by Taubes and Perutz, we study the geometry of {\it near-symplectic} $4$--manifolds, i.e., manifolds equipped with a closed $2$--form which is symplectic outside a union of embedded $1$--dimensional submanifolds, and {\it broken Lefschetz fibrations} on them \cite{adk,gkirby}. We present a set of four moves which allow us to pass from any given broken fibration to any other which is deformation equivalent to it. Moreover, we study the change of the near-symplectic geometry under each of these moves. The arguments rely on the introduction of a more general class of maps, which we call {\it wrinkled fibrations} and which allow us to rely on classical singularity theory. Finally, we illustrate these constructions by showing how one can merge components of the zero-set of the near-symplectic form. We also disprove a conjecture of Gay and Kirby by showing that any achiral broken Lefschetz fibration can be turned into a broken Lefschetz fibration by applying a sequence of our moves.

\end{abstract}

\maketitle
\section{Introduction}
\label{section1}

\subsection{Near-symplectic manifolds}
Let $X$ be a smooth, oriented $4$--manifold. Then a closed $2$--form $\omega$ is called {\it near-symplectic} if $\omega^2\geq 0$ and there is a metric $g$ such that $\omega$ is self-dual harmonic and transverse to the $0$--section of $\Lambda^+$. Equivalently, without referring to any metric, one could define a closed $2$--form $\omega$ to be near-symplectic if for any point $x \in X$ either $\omega^2_x > 0$,  or $\omega_x=0$, and the intrinsic gradient $(\nabla\omega)_x\colon T_x X \to \Lambda^2 T_x^*X$ has maximal rank, which is $3$.

The zero-set $Z$ of such a $2$--form is a $1$--dimensional submanifold of $X$. If $X$ is compact and $b_2^+(X)>0$ then Hodge theory gives a near-symplectic form $\omega$ on $X$. Clearly, in this case $Z$ is just a collection of disjoint circles. Furthermore, by deforming $\omega$, one can show that on any near-symplectic manifold, one can reduce the number of circles to $1$, this was proved in \cite{perutz}. We give a new proof of this result in Theorem \ref{6.1} as an application of the techniques developed in this paper. Of course, the last circle cannot be removed unless the underlying manifold is symplectic.

Interesting topological information about $X$ is captured by the natural decomposition of the normal bundle of these circles, provided by the near-symplectic form. More precisely, transversality of $\omega$ implies that $\nabla\omega \colon N_Z \to \Lambda^+ X $ is an isomorphism, where $N_Z$ is the normal bundle to the zero-set of $\omega$. This enables us to orient the zero-set $Z$. Now consider the quadratic form $N_Z \to \f{R}$, $v \to \langle\iota(z)\nabla_v \omega, v \rangle $, where $z$ is a non-vanishing oriented vector field on $Z$. As $dw=0$, this quadratic form is symmetric and has trace zero. It follows that, it has three real eigenvalues everywhere, where two are positive and one is negative. Then $N_Z=L^+{\oplus} L^-$, where $L^{\pm}$ are the positive and negative eigen-subbundles respectively. 
In particular, this allows us to divide the zero-set into two pieces, the even circles where the line bundle $L^-$ is orientable, and the odd circles where $L^-$ is not orientable. This definition is motivated by the following result of Gompf that the number of even circles is equal to $1-b_1 + b_2^+ $ modulo $2$  \cite{perutz} . In particular, observe that if there is only one zero circle which is even, the manifold $X$ cannot be symplectic.

In this paper, we will be interested in local deformations of near-symplectic forms on a $4$--manifold. An important such deformation is provided by the Luttinger--Simpson model given on $D^4 \subset \f{R}^4$ where the birth (or death) of a circle can be observed explicitly \cite{perutz}:\\
\begin{eqnarray*} 
 \omega_s & = & 3 \epsilon (x^2+t^2-s) (dt \wedge dx + dy \wedge dz ) + 6 \epsilon y (tdt \wedge dz + x dx \wedge dz) \\ 
  & - & 2z (dx \wedge dy + dt \wedge dz ) + 2y (dt \wedge dy + dz \wedge dx ) \\
\end{eqnarray*}
for $\epsilon \leq \frac{1}{6}$.

We will see that this is not the only type of deformation of near-symplectic forms. One of the goals of this paper is to identify such deformations and interpret them in terms of the singular fibrations associated to them. 


\subsection{Wrinkled fibrations}

A \emph{broken fibration} on a closed $4$--manifold $X$ is a smooth map to a closed surface with singular set $A \sqcup B$, where $A$ is a finite set of singularities of Lefschetz type where around a point in $A$ the fibration is locally modeled in oriented charts by the complex map $(w,z) \to w^2+z^2 $, and B is a $1$--dimensional submanifold along which the singularity of the fibration is locally modeled by the real map $(t,x,y,z) \to (t, x^2+y^2 - z^2) $, $B$ corresponding to $t=0$. We remark here that we do not require the broken fibrations to be embeddings when restricted to their critical point set. In particular, this means that the critical value set may include double points.

There have been two different approaches to constructions of broken fibrations on $4$--manifolds. The first approach is by Auroux, Donaldson and Katzarkov \cite{adk} based on approximately holomorphic techniques, generalizing the construction of Lefschetz pencils on symplectic manifolds. The more recent approach is due to Gay and Kirby \cite{gkirby}, where the fibration structure is constructed explicitly in two pieces in the form of open books, and then Eliashberg's classification of overtwisted contact structures as well as Giroux's theorem of stabilization of open books are invoked to glue these two pieces together to form an achiral broken fibration. Achiral here refers to the existence of finitely many Lefschetz type singularities with the opposite orientation on the domain, namely the singularity is modeled by the complex map $(w,z) \to \bar{w}^2+z^2 $. 

There is a correspondence between broken fibrations and near-symplectic manifolds up to blow-up, in analogy with the correspondence between Lefschetz fibrations and symplectic manifolds up to blow-up. More precisely, given a broken fibration on a $4$--manifold $X$ with the property that there is a class $h\in H^2(X)$ such that $h(F)>0$ for every component $F$ of every fibre, it is possible to find a near-symplectic form on $X$ such that the regular fibres are symplectic and the zero-set of the near-symplectic form is the same as the $1$--dimensional critical point set of the broken fibration.  This is an adaptation due to Auroux, Donaldson and Katzarkov of Gompf's generalization of Thurston's argument used in finding a symplectic form on a Lefschetz fibration. Conversely, in \cite{adk}, it is proven that on every $4$--manifold with $b_2^+ (X ) > 0 $ (recall that this is equivalent to $X$ being near-symplectic), there exists a broken fibration if we blow up enough. One of the questions of interest that remains to be answered is to determine how unique this broken fibration is. In particular, we would like to find a set of moves on broken fibrations relating two different broken fibrations on a given $4$--manifold. One of the main themes of this paper is the discussion of a set of moves which allows one to pass from one broken fibration structure to another.

In this paper, we will consider a slightly more general type of fibration, where we will allow cuspidal singularities on the critical value set of the fibration. These type of fibrations occur naturally when one considers deformations of the broken fibrations. We will also discuss a local modification of a cuspidal singularity (without changing the diffeomorphism type of the underlying manifold structure) in order to get a broken fibration. Therefore, one can first deform a broken fibration to obtain a wrinkled fibration, then apply certain moves to this wrinkled fibration, and finally modify the wrinkled fibration in a neighborhood of cuspidal singularities to get a genuine broken fibration. In this way, one obtains a set of moves on broken fibrations on a given $4$--manifold.

Let $X$ be a closed $4$--manifold, and $\Sigma$ be a $2$--dimensional surface. We say that a map $f\colon X\to \Sigma$ has a \emph{cusp singularity} at a point $p\in X$, if around $p$, $f$ is locally modeled in oriented charts by the map $(t,x,y,z) \to (t,x^3-3xt+y^2-z^2)$. This is what is known as the Whitney tuck mapping, the critical point set is a smooth arc, $\{x^2=t, y=0,z=0\}$, whereas the critical value set is a cusp, namely it is given by $C = \{ (t,s) \colon 4t^3 = s^2 \} $. This is the generic model for a family of functions $\{f_t\}$, which are Morse except for finitely many values of $t$ \cite{arnold}. The signs of the terms $y^2$ and $z^2$ are chosen so that the functions $f_t$  have only index $1$ or $2$ critical points.
More precisely, if $f\colon\f{R}^3\to\f{R}$ is a Morse function with only index $1$ or $2$ critical points, then $F\colon\f{R}^4\to \f{R}^2 $ given by $(t,x,y,z)\to (t, f(x,y,z))$ is a broken fibration with critical set in correspondence with the critical points of $f$. Notice that the functions $f_t(x,y,z)=x^3-3xt+y^2-z^2$ are Morse except at $t=0$, where a birth of critical points occur.
\ \\
\begin{definition} A wrinkled fibration on a closed $4$--manifold $X$ is a smooth map $f$ to a closed surface which is a broken fibration when restricted to $X \backslash C$, where $C$ is a finite set such that around each point in $C$, $f$ has cusp singularities. We say that a fibration is {\it purely wrinkled} if it has no isolated Lefschetz-type singularities.
\end{definition}

It might be more appropriate to call these fibrations ``broken fibrations with cusps'', to avoid confusion with the terminology introduced by Eliashberg-Mishachev \cite{yasha}. The reason for our choice of terminology is that wrinkled fibrations can typically be obtained from broken Lefschetz fibrations by applying wrinkling moves  (see move $4$ in Section \ref{section3}) which eliminates a Lefschetz type singularity and introduces a wrinkled fibration structure. Conversely, as mentioned above, it is possible to locally modify a wrinkled fibration by smoothing out the cusp singularity at the expense of introducing a Lefschetz type singularity and hence get a broken fibration. 

\begin{theorem}
\
\begin{enumerate}
\item[a)]	
	Every wrinkled fibration is homotopic to a broken fibration by a homotopy supported near cusp singularities. 
\item[b)]
	Every broken fibration is homotopic to a purely wrinkled fibration by a homotopy supported near Lefschetz singularities. 
\end{enumerate}	
\end{theorem}

The first part of this paper is concentrated on a set of moves on wrinkled fibrations and corresponding moves on broken fibrations. All of these moves keep the diffeomorphism type of the total space unchanged. We remark here that as will be explained below these moves occur as deformations of wrinkled fibrations and not as deformations of broken fibrations. To be more precise, by a {\it deformation} of wrinkled fibrations we mean a one-parameter family of maps which is a wrinkled fibration for all but finitely many values of the parameter. In fact, as we will see, an infinitesimal deformation of a broken fibration gives a wrinkled fibration whereas the wrinkled fibrations are stable under infinitesimal deformations. This is indeed the main reason for extending the definition of the broken fibrations to wrinkled fibrations. 

In the second part, using techniques from singularity theory, we prove that our list of moves is complete in the sense that any generic infinitesimal deformation of a wrinkled fibration which does not have any Lefschetz type singularity is given by one of the moves that we exhibited in Section \ref{section3}. Furthermore, as we will see in Section \ref{section3}, it is always possible to deform a wrinkled fibration infinitesimally so that the Lefschetz type singularities are eliminated. \ \\

\begin{theorem}
\
\begin{enumerate}
\item[a)] Any one-parameter family deformation of a purely wrinkled fibration is homotopic rel endpoints to one which realizes a sequence of births, merges, flips, their inverses and isotopies staying within the class of purely wrinkled fibrations.
\item[b)] Given two broken fibrations, suppose that after perturbing them to purely wrinkled fibrations, the resulting fibrations are deformation equivalent. Then one can get from one broken fibration to the other one by a sequence of birth, merging, flipping and wrinkling moves, their inverses and isotopies staying within the class of broken fibrations. 
\end{enumerate}
\end{theorem}

As in the case of broken fibrations, one can define a wrinkled pencil on $X$ to be a wrinkled fibration $f\colon X\backslash P \to \Sigma$, where $P$ is a finite set and around a point in $P$, the fibration is locally modeled in oriented charts by the complex map $(w,z) \to w/z$. Note that, after blowing up $X$ at the points $P$, one can get a wrinkled fibration. It is possible to construct a natural near-symplectic form that is ``adapted'' to a given wrinkled pencil. The key property of this form is that it should restrict to a symplectic form on the smooth fibres of the given wrinkled fibration. Therefore, we can equip every wrinkled pencil with a well-defined deformation class of near-symplectic forms, it is natural thus to study what happens to this class after each move that was described on the previous paragraph. This will be discussed Section \ref{section5} of the paper.
  
In Section \ref{section6}, we give a number of applications of our moves on broken fibrations. Notably, by considering the mirror image of the wrinkling move, we prove that we can turn an achiral Lefschetz singularity into a wrinkled map and then into a broken fibration, without losing equatoriality of the round handles. This provides the following simplification of the result of Gay and Kirby in \cite{gkirby}:  

\begin{theorem} Let $X$ be an arbitrary closed $4$--manifold and let $F$ be a closed surface in $X$ with $F\cdot F = 0$. Then there exists a broken Lefschetz fibration from $X$ to $S^2$ with embedded singular locus, and having $F$ as a fibre. Furthermore, one can arrange so that the singular set on the base consists of circles parallel to the equator with the genera of the fibres in increasing order from one pole to the other.

\end{theorem}

We remark that this disproves the conjecture $1.2$ of Gay and Kirby in \cite{gkirby} about the essentialness of including achiral Lefschetz singularities for broken fibrations on arbitrary closed $4$--manifolds. 

After the first writing of this paper, an earlier result of a similar nature, but allowing the set of critical values of the fibration to be immersed rather than embedded, has been obtained by Baykur in \cite{baykur}. Namely, Baykur proved an existence theorem for broken fibrations with immersed critical value set by combining the following two ingredients: (1) a result of Saeki \cite{Saeki} which says that any continuous map from a closed $4$-manifold $X \to S^2$ is homotopic to a stable map without definite folds, i.e. in our terminology, a purely wrinkled fibration with immersed critical value set, (2) the cusp modification described in Section $2$ of this paper. 

Another recent development that took place after the writing of this paper is worth mentioning here: Akbulut and Karakurt \cite{AkKara} came up with a new proof of the existence theorem stated above by refining the construction of Gay and Kirby. The difference between Akbulut and Karakurt's result and ours is that they directly construct a broken fibration on any $4$-manifold, whereas we describe a way to modify achiral Lefschetz singularities into broken and Lefschetz singularities.

Finally, here we would like to discuss our main motivation for studying the particular structure of broken fibrations and their deformations, the wrinkled fibrations.

\subsection{Seiberg--Witten invariants and Lagrangian matching invariants}

In \cite{DonaldsonSmith}, Donaldson and Smith define an invariant of a symplectic manifold $X$ by counting holomorphic sections of a relative Hilbert scheme that is constructed from a Lefschetz fibration on a blow-up of $X$. More precisely, by Donaldson's celebrated theorem, there exists a Lefschetz fibration $f\colon X'\to S^2$, where $X'$ is some blow-up of $X$. Then, for any natural number $r$, Donaldson and Smith give a construction of a relative Hilbert scheme $F\colon X_r(f)\to S^2$, where the fibre over a regular value $p$ of $f$ is the symmetric product $\Sigma^r(f^{-1}(p))$. In fact, $X_r(f)$ is a resolution of singularities for the relative symmetric product, which is the fibration obtained by taking the $r^{th}$ symmetric product of each fibre. They then define their standard surface count, which is some Gromov invariant counting pseudoholomorphic sections of $X_r(f)$. Usher, in \cite{usher}, proves that this invariant is the Gromov invariant of the underlying symplectic $4$--manifold $X$. Finally, we know that this is in turn equal to the Seiberg-Witten invariant of $X$ by the seminal work of Taubes \cite{taubes}. Therefore, one obtains a geometric formulation of the Seiberg-Witten invariant for a symplectic manifold $X$ on a Lefschetz fibration structure associated to $X$, which also shows in particular that this invariant is independent of the Lefschetz fibration structure. 

A similar but technically not so straightforward generalization of this method of getting an invariant from a Lefschetz fibration is described in \cite{matching} for the case of broken fibrations, thus giving an invariant for all smooth $4$--manifolds with $b_2^+(X) > 0$. These are called Lagrangian matching invariants. Here we give a quick sketch of the definition of these invariants.

Suppose $X$ is a near-symplectic manifold with only one zero circle $Z$, and $f\colon X\to S^2$ is a broken fibration with one circle of singularity along the equator of $S^2$. Take out a thin annulus neighborhood of the equator and write $N$ and $S$ for the closed discs that contain the north pole and the south pole respectively. Let $X^N=f^{-1}(N)$ and $X^S=f^{-1}(S)$, suppose the fibre genus of $X^N$ is $g$ and the fibre genus of $X^S$ is $g-1$. Consider the relative Hilbert schemes Hilb$^r_N(X^N)$ and Hilb$^{r-1}_S(X^S)$. These are symplectic manifolds with boundaries $Y^N_r = \Sigma^r_{S^1}(\del X^N)$ and $Y^S_{r-1} = \Sigma^{r-1}_{S^1}(\del X^S)$, respectively. 

Perutz then constructs a sub-fibre bundle $\mathcal{Q}$ of the fibre product $Y^N_r \times_{S^1} Y^S_r \to S^1$ which constitutes the Lagrangian boundary conditions for the pairs of pseudoholomorphic sections of Hilb$^r_N(X^N)$ and Hilb$^{r-1}_S(X^S)$ in the following sense : One defines $\mathcal{L}_{X,f}$ to be a Gromov invariant for pairs $(u_N, u_S)$ of pseudoholomorphic sections of Hilb$^r_N(X^N)$ and Hilb$^{r-1}_S(X^S)$ such that the boundary values $(u_N|_{\del N}, u_S|_{\del S})$ lie in $\mathcal{Q}$.

Now, the big conjecture in this field is of course the conjecture that Lagrangian matching invariants equal the Seiberg-Witten invariants. This has been verified by Perutz, \cite{matching}, in several cases, notably in the case of symplectic manifolds as mentioned above, and when the underlying manifold is of the form $S^1 \times M^3$, for any $M$ which is a $\f{Z}$--homology--$(S^1\times S^2)$ and for connected sums.

An important problem to be explored is that the Lagrangian matching invariant is not yet known to be an invariant of the given $4$--manifold. In other words, it is an invariant of the near-symplectic manifold together with a given broken fibration structure. Our next task in this field will be to show that the Lagrangian matching invariant stays an invariant under the set of moves that we describe in this paper. We believe that our set of moves will be enough to pass from a given broken fibration structure on a manifold to any other broken fibration structure on the same manifold under suitable hypotheses on the homotopy type of the fibration map. We have strong evidence for this since, as was mentioned above, the set of moves that we discuss in this paper are sufficient to pass from a given broken fibration to any one-parameter deformation of it. These two hypotheses would imply that the Lagrangian matching invariant is really an invariant of the underlying manifold. We believe that these steps will play an important role in proving the big conjecture mentioned above.

\textbf{Acknowledgements : } The author is grateful to his Ph.D.\ supervisor Denis Auroux for valuable suggestions and comments. He is also indebted to Tim Perutz for encouraging and helpful remarks. Thanks also to Karl Luttinger for sharing a copy of the unpublished manuscript \cite{luttinger} and to Robion Kirby for his interest in this work. This work was partially supported by an NSF grant (DMS-0600148).

\section{A local modification on wrinkled fibrations}
\label{section2}
Recall from the introduction that a cusp singularity is given locally by the map $F \colon \f{R}^4 \to \f{R}^2 $ given by: \[ (t,x,y,z)\to (t, x^3 - 3xt + y^2 - z^2) \]

The critical point set is a smooth arc, $\{x^2=t, y=0,z=0\}$, whereas the critical value set is a cusp, namely it is given by $C = \{ (t,s) \colon 4t^3 = s^2 \} $.

The idea is to modify a neighborhood of the singular point of the cusp with an allowed model for broken fibrations without changing the topology. We will do this by surgering out a neighborhood of the cuspidal singularity and gluing back in a neighborhood of an arc together with a Lefschetz type singular point as shown in Figure $1$. The issue is to make sure that the fibration structures match outside the neighborhood.

\begin{figure}[h]
\begin{center}
\includegraphics{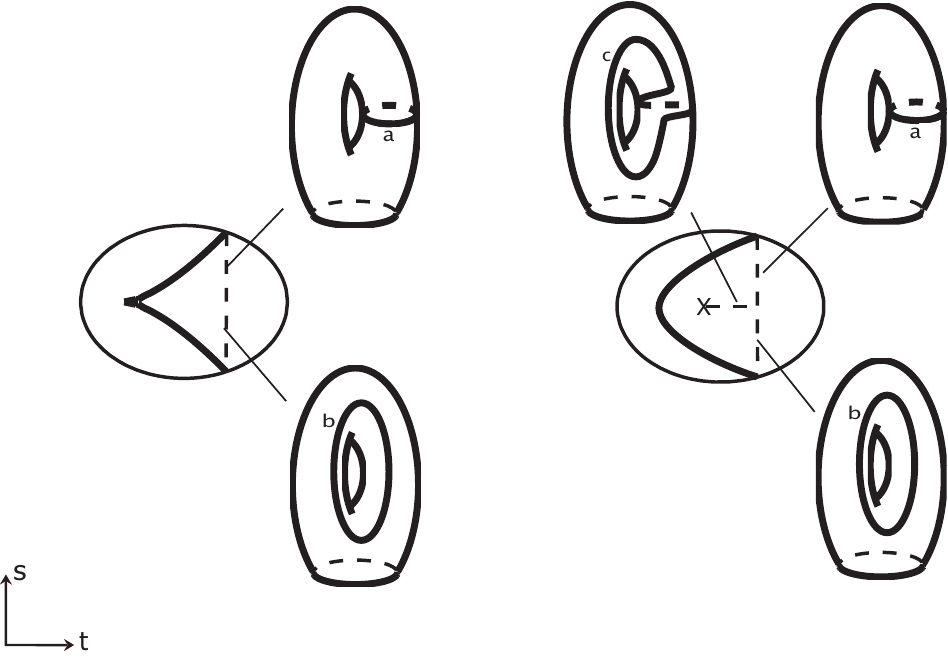}
\end{center}
\caption{Local modification}
\end{figure}

Restricting to a neighborhood of the origin, we get a map $F \colon D^4 \to D^2 $ and $C$ divides the image into two regions, where the fibres above the ``interior region'' are punctured tori, whereas the fibres above the ``exterior region'' are discs, as shown in Figure $1$. Furthermore, looking above the line $\{ t = \frac{1}{2} \}$, one sees that as the parameter $s$ converges to $C$ from below, $s\to \frac{1}{\sqrt{2}}$,  one of the generating loops of the homology of the torus collapses to a point , and as $s$ converges to $C$ from above, $s\to -\frac{1}{\sqrt{2}}$ the other generator collapses to a point. This is evident from the fact that $f_{1/2}(x,y,z) = x^3 - \frac{3}{2} x + y^2 - z^2 $ restricted to the preimage of $\{ t = \frac{1}{2} \}$ is a Morse function on $D^3$ with $2$ critical points of indices $1$ and $2$ which cancel each other. 
 
Now consider the $D^2$--valued broken fibration structure described on the right of Figure $1$. Let us denote this fibration by $p \colon X \to D^2$. This fibration is cooked up so that it matches above a neighborhood of the boundary of $D^2$ with the fibre structure of the map $F$. On the other hand, by introducing a Lefschetz type singularity, we are able to have a broken fibration structure, where the vanishing cycles are described on the right of the Figure $1$.
In order to perform a local surgery to pass from the map $F$ to the described broken fibration $p$, it remains to show that the total space $X$ is diffeomorphic to $D^4$. This will be accomplished by giving a handle decomposition of $X$, and showing that it is in fact obtained by attaching one $1$--handle and one $2$--handle to $D^4$, in such a way that they can be cancelled.

Let us now describe $X$ explicitly. Denote the standard loops generating homology of a regular fibre by $a$ and $b$. As shown in Figure $1$, restricting to the line $\{ t = \frac{1}{2} \}$, as $s$ approaches to $C$ from below, $a$ collapses to a point and as $s$ approaches to $C$ from above, $b$ collapses. (This is to be consistent with the fibre structure of $F$.) Now the monodromy around the Lefschetz type singularity must be the Dehn twist along $c=a-b$, denoted by $\tau_{a-b}$ so that $\tau_{a-b}(a) = b$, where we oriented $a$ and $b$ so that $a\cdot b = -1$. Therefore, restricting to the line $s=0$, as $t$ approaches to the singularity $a-b$ collapses to a point. (Here by $c=a-b$ we really denote an embedded loop $c$ which is equal to $[a-b]$ as a homology class.) We remark here that, just as in Lefschetz fibrations, a diagram indicating the fibre structure and vanishing cycles along relevant paths is enough to determine a broken (or wrinkled) fibration uniquely on a disc. We now have an explicit understanding of the various vanishing cycles for X. Next we proceed to describe the corresponding handle diagram. We first restrict to the preimage of the region shown in Figure $2$. This is clearly diffeomorphic to the total space $X$. Now divide this region into $3$ parts as shown in Figure $2$. The preimage of region $0$ is just $D^2\times D^2 = D^4$. We claim that the preimage of regions $0$ and $1$ together is $D^4\cup$ $1$--handle, and the preimage of all three regions is $D^4 \cup 1$--handle $ \cup$ $2$--handle in such a way that the attaching sphere of the $2$--handle intersects the belt sphere of the $1$--handle transversely at a single point, so that these two handles can be cancelled. 

\begin{figure}[!ht]
\begin{center}
\includegraphics{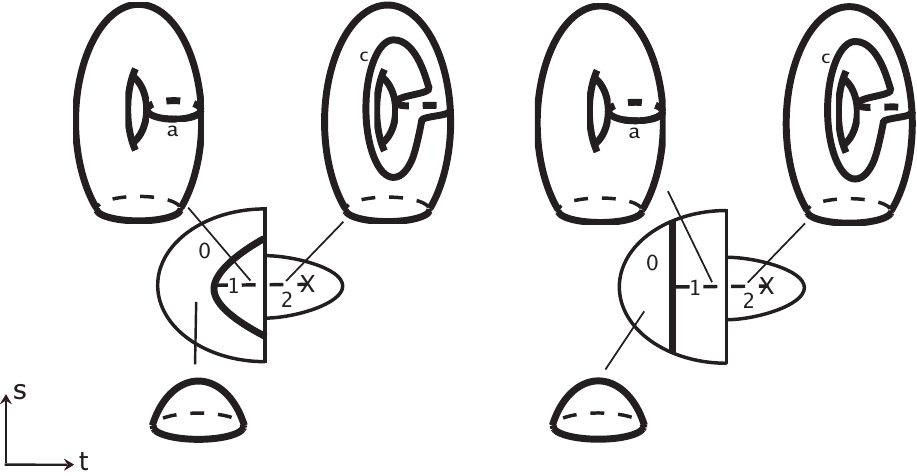}
\end{center}
\caption{Handle decomposition of the total space}
\end{figure}

In this picture, it is more convenient to fix the reference fibre above a point which lies between regions $1$ and $2$ as shown in the Figure $2$. Just for simplicity, we can choose an identification of this reference fibre with the previous choice using the parallel transport along a simple arc above the Lefschetz singularity so that the vanishing cycles in this new reference fibre are given as shown. Finally, observe that we can isotope the base so that the $1$--dimensional singular set is straightened to a line. 

Next, we are in a position to see the handle decomposition very explicitly. In fact, the preimage of the regions $0$ and $1$ can be thought as ( $D^3\ \cup$ $1$--handle $ ) \times D^1 $, where the $D^1$ is the $s$ direction. The belt circle of this $3$--dimensional $1$--handle corresponds to the vanishing cycle $a$ on a regular fibre above the region $1$, to be precise, fix the regular fibre $F$ above a point $p$ in region $1$, say $p$ lies on the $s=0$ line. Now ( $D^3\ \cup$ $1$--handle $ ) \times D^1 = D^4\ \cup$ $1$--handle where the belt sphere of this latter $4$--dimensional $1$--handle intersects $F$ at $a$. Now, by construction starting from $F$ as one approaches to Lefschetz singularity the loop $c$ collapses to a point. It is a standard fact of Lefschetz fibrations that gluing the preimage of region $2$ corresponds to a $2$--dimensional handle attachment with attaching circle being the loop $c$ on $F$ \cite{gompfstip}. (In fact, if one considers the local model $(z,w)\to z^2+w^2$, then $Re(z^2+w^2)$ is a Morse function with one critical point of index $2$ at the origin.) Therefore, we conclude that $X = D^4\ \cup$ $1$--handle $\cup$  $2$--handle with belt sphere of the $1$--handle intersects the attaching circle of the $2$--handle transversely at exactly one point, and this intersection point is precisely  the intersection of the loop $a$ and the loop $c$ on $F$. Finally, applying the cancellation theorem of handle attachments, we conclude that $X = D^4$ as required.

It is of interest to note that one could as well replace a cusp singularity with a broken arc singularity and an achiral Lefschetz singularity, where vanishing cycles for the cusp are given by $a$ and $b$ as before, and the vanishing cycle for the achiral Lefschetz singularity is given by $c= a+b$ (since one must now have $\tau_c^{-1}(a) = -b $). The difference between a Lefschetz singularity and an achiral Lefschetz singularity with the same vanishing cycle is in the framing of the corresponding $2$--handle attachment. Namely, a Lefschetz singularity corresponds to $-1$--framing with respect to the fibre framing whereas an achiral Lefschetz singularity corresponds to $+1$--framing. The cancellation theorem of handle attachments does not see the framings, therefore the proof is verbatim.

We remark here that the local modification described in this section is not given as a deformation, in the sense that we have not explained how to give a one-parameter family of wrinkled fibrations which starts from the fibration depicted on the left side of Figure $1$ and ends at the fibration given on the right side of Figure $1$. We will actually give such a family in the next section, which will in fact give yet another way of proving the validity of the above move. However, we chose to present the above proof first, as it is considerably simpler and in fact this enabled the author to discover more complicated modifications described in Section \ref{section3}, which come equipped with deformations. Afterwards, we were able to recover the above modification as a composition of these deformations. 
\section{A set of deformations on wrinkled fibrations}
\label{section3}
In this section, we describe a set of moves on wrinkled fibrations. We first give three such moves which are deformations of wrinkled fibrations and the corresponding deformations which end up being broken fibrations are obtained by applying the modification described in Section \ref{section2}, which as was mentioned there, is indeed a deformation. Note that this was not proved in the previous section. This will be accomplished after we describe the last move which enables us to turn a Lefschetz singularity into three cusp singularities.
\paragraph{Move 1 (Birth) :} Consider the wrinkling map $F_s \colon \f{R}\times \f{R}^3 \to \f{R}^2$, as defined in \cite{yasha} :
\[ (t,x,y,z) \to (t, x^3+3(t^2-s)x+y^2-z^2)\]

For $s<0$, this is a genuine fibre bundle, i.e., there is no singularity. At $s=0$, the only singularity is at the origin. This is a degenerate map, which is not an allowed singularity for a wrinkled fibration. For $s>0$, the critical point set of $F_s$ is the circle $\{x^2+t^2=s, y=z=0 \}$, whereas the critical value set $C_s$ is a \emph{wrinkle} shown on the left of Figure $3$. This is clearly a wrinkled map. Therefore, we have a deformation of wrinkled maps, the only subtle change being at $s=0$, where birth of the wrinkle happens. 

\begin{figure}[!ht]
\begin{center}
\includegraphics{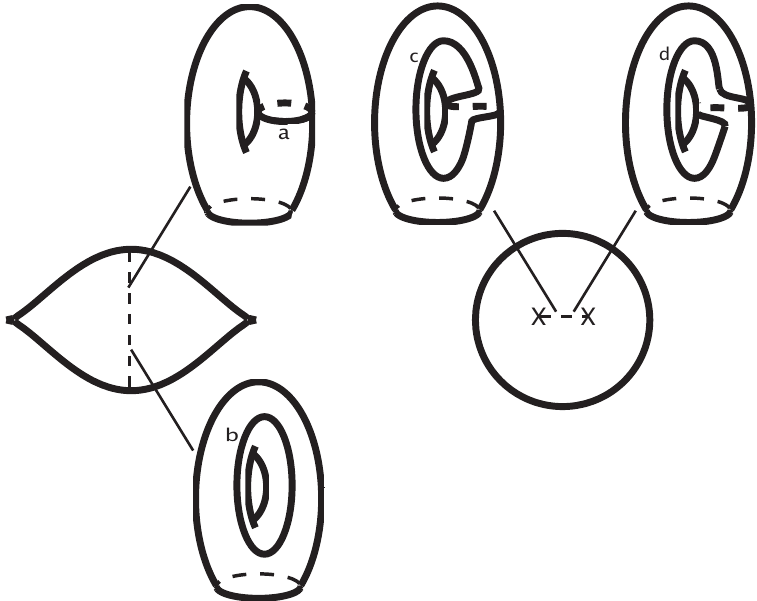}
\end{center}
\caption{Creation of a circle singularity along with two point singularities}
\end{figure}

Now, fix $s=1$. Considering the wrinkle as obtained from gluing two cusps together, we can apply the local modification of Section \ref{section2} to obtain a broken fibration on $\f{R}^4$ with singular set a circle together with two point singularities as shown on the right of Figure $3$. Note that one has to check that the configuration of the vanishing cycles matches the model in Section \ref{section2}. Conveniently, we can check this on the vertical line ${t=0}$. Then the map becomes $(0,x,y,z) \to (0,x^3-3x+y^2-z^2)$, and this is the same map that was used in Section \ref{section2}, therefore the configuration of the vanishing cycles matches the model in Section \ref{section2}. Namely, on ${t=0}$, the two vanishing cycles obtained from approaching to $C$ from below and from above starting from the origin, intersect transversely at a point.

Thus, given a broken fibration on any $4$--manifold, we can restrict the fibration to a $D^2$ on the base where the fibration is regular, and also restrict the fibres to obtain $D^2\times D^2$. Then, apply the move just described to obtain a new fibration, where the singular set is changed by an addition of a circle and two points. Furthermore, the fibre genus above the points in the interior of this new singular circle increases by $1$. 

We remark that this move on broken fibrations was first observed by Perutz in proposition $1.4$ of \cite{matching}, where he proves that the total space of the closed fibre case of the fibration on the right of Figure $3$ is diffeomorphic to $S^2 \times S^2$. Here, we were able to divide this move into two pieces by allowing cusp singularities, which indicates that the local move of Section \ref{section2} is a more basic move. 

\paragraph{Move 2 (Merging) :} Let us now describe another move which corresponds to merging two singular circles to obtain one circle together with two Lefschetz type singularities. We begin with the local picture described on the left side of Figure $4$. The lines which separate regions on the base indicate the critical value set. The vanishing cycles obtained from moving towards the upper line and moving towards the lower line are assumed to intersect transversely at a singular point. The standard model for such a broken fibration $F\colon\f{R}^4\to\f{R}^2$ is given by the map $(t,x,y,z) \to (t, x^3-3x+y^2-z^2)$. The critical value set of this map is given by two horizontal lines, and the configuration of the vanishing cycles is as described. Now consider the map $F_s \colon \f{R}\times \f{R}^3 \to \f{R}$ defined by:
\[ (t,x,y,z) \to (t, x^3+3(s-t^2)x+y^2-z^2)\]

Then for $s<0$, $F_s$ is isotopic to $F$, with the critical value set 
being $C = \{ (t,u) \colon 4(t^2-s)^3 = u^2 \} $. For $s<0$, $C$ consists of two simple curves and is isotopic to the left side of Figure $4$. At $s=0$, as before, we have a more degenerate map. This is where a subtle change in the fibration structure occurs. For $s>0$, we get a wrinkled map with critical value set, including two cusp singularities, isotopic to the model depicted in the middle part of Figure $4$. Note that the picture on Figure $4$ is drawn so that the maps are equal outside of a neighborhood of the origin, to ensure that when restricted to $D^4$, the maps agree on a neighborhood of the boundary. Finally, we apply the local modification model from Section \ref{section2} to each cuspidal singularity to get a new broken fibration. Therefore, we obtain a move of broken fibrations, namely whenever one has the configuration described on the left side of Figure $4$, one can surger out a $D^4$ and glue in the right side of Figure $4$ to obtain a new configuration.

\begin{figure}[!ht]
\begin{center}
\includegraphics{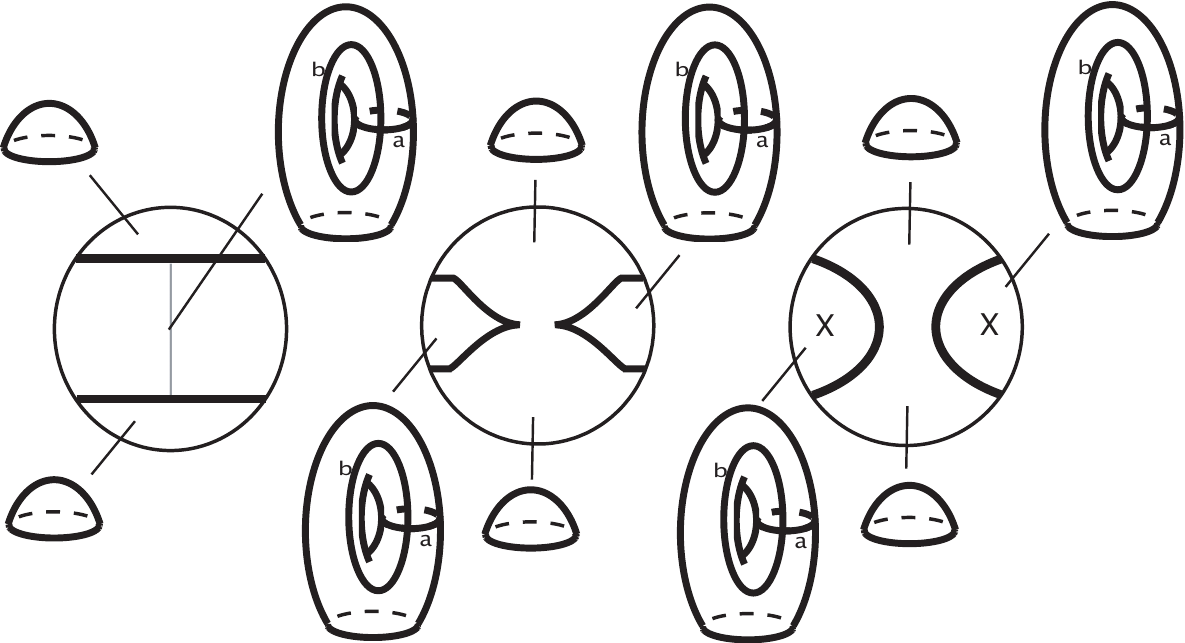}
\end{center}
\caption{Merging singular circles}
\end{figure}

We remark here that to apply a merging move, one needs a configuration as in the left side of Figure $4$, in particular it is necessary that the vanishing cycles intersect transversely at a unique point. On the other hand, to apply an inverse merging move the following two conditions are necessary. Referring to the right part of Figure $4$, one needs to make sure that, fixing a reference fibre halfway along a path connecting the Lefschetz singularity and the broken singularity on the left, the vanishing cycles for the Lefschetz singularity and the broken singularities should intersect transversely at a point. Exactly the same configuration is required on the right side of the fibration. However, we would like to point out that there is no compatibility condition required for the two sides as long as the fibres in the middle region are connected. Namely, to give an embedding of the fibration depicted on the right side of Figure $4$ into a fibration that has the same base and whose vanishing cycles satisfy the condition described above, one divides the base into three pieces: a piece on the left that includes the Lefschetz singularity and the broken singularity, a middle piece which is a smooth fibration, and a piece on the right which includes the Lefschetz singularity and the broken singularity. Since the vanishing cycles are as prescribed above, it is easy to construct a fibrewise embedding of the total spaces of the pieces on the left and on the right. Namely, given two simple closed curves intersecting transversely at one point on a fibre $F$, it is always possible to find a diffeomorphism of $F$ such that those two curves are standardized, in the sense that they sit in the standard way as part of an embedding of a punctured torus to $F$. Finally in order to give an embedding of the total space of the middle piece, one needs to give a fibrewise embedding of the disc fibration $D^2 \times D^2$ such that, if we consider the base $D^2$ as $[0,1]^2$, the embedding is already prescribed above $\{0,1\}\times [0,1]$. But now, it is easy to extend this to a fibrewise embedding of $D^2 \times D^2$ by just flowing the fibers above $\{0\}\times[0,1]$ to fibres above $\{1\}\times [0,1]$ since the set of embeddings of $D^2$ to a fibre $F$ is clearly connected  provided that $F$ is connected.

\paragraph{Move 3 (Flipping) :} This move is originally due to Auroux. The observation was that for a given near-symplectic manifold $(X,\omega)$, if one considers possible broken fibrations adapted to $(X,\omega)$, the rotation number of the image of a given component of the zero-set of $\omega$ is not fixed a priori. If one considers a one-parameter family of deformations of broken fibrations, one can possibly get a flip through a real cusp. However, here we discuss this move in an alternative way to the original approach, using the local modification discussed in Section \ref{section2}. Consider the map $F_s \colon \f{R}\times \f{R}^3 \to \f{R}^2$ given by:
\[ (t,x,y,z) \to (t, x^4-x^2 s+xt+y^2-z^2)\]
Then for $s<0$, the critical value set consists of a simple curve and $F_s$ is isotopic to the map described on the left side of Figure $5$. At $s=0$, we have a higher order singularity and as before this is where a subtle change in the fibration structure occurs. For $s>0$, we get a wrinkled map with critical value set, including two cusp singularities, isotopic to the model depicted in the middle part of Figure $5$. This map still induces an immersion on the critical point set away from the cusp singularities, however now we have a double point as shown in Figure $5$.

\begin{figure}[!ht]
\begin{center}
\includegraphics{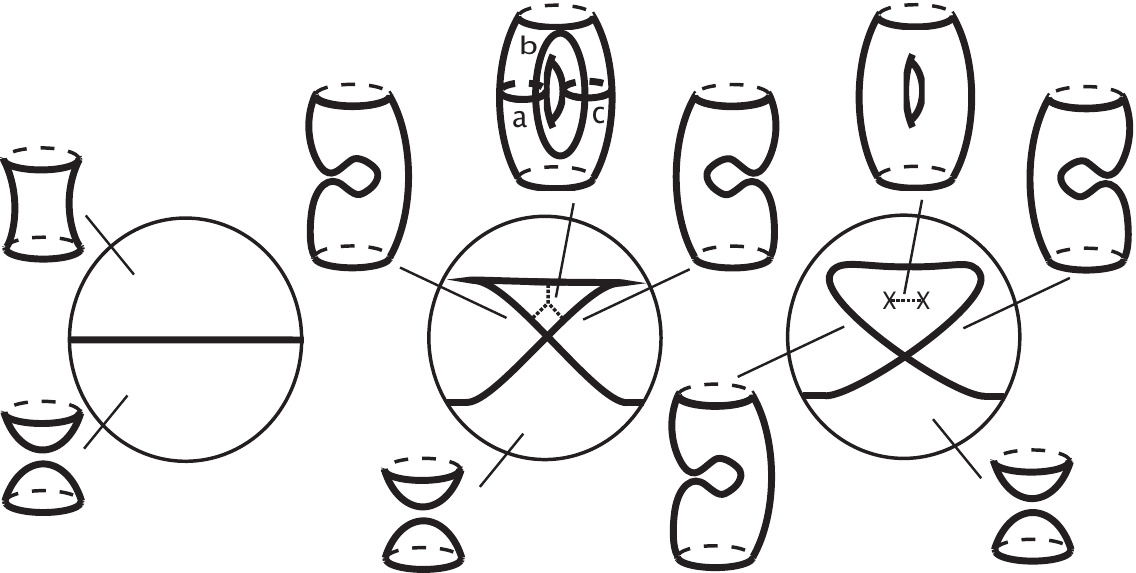}
\end{center}
\caption{Flipping}
\end{figure}

One can fix a reference fibre in the interior region (the high-genus region) as in the middle portion of Figure $5$ so that the vanishing cycles for the three paths drawn are the given loops $a$, $b$, $c$. Indeed, we know from the local model of a cusp singularity that the vanishing cycles corresponding to each branch of a cusp intersect transversely once. Therefore, the vanishing cycle for the path going up intersects both the vanishing cycle for the lower left path and the vanishing cycle for the lower right path transversely at a point. Furthermore, we know that the two latter vanishing cycles are disjoint since the critical point set in the total space is embedded, and they cannot be homotopic, since otherwise the fibres above the bottom region would have a sphere component. Now, once these intersection properties are understood, it is easy to see that there is a diffeomorphism of the twice punctured torus that sends any configuration of three simple closed curves satisfying the above properties to $a$, $b$ and $c$.

On the right side of Figure $5$, it follows from monodromy considerations (recall that the monodromy around a Lefschetz singularity is the Dehn twist along the vanishing cycle) as in Section \ref{section2} that the vanishing cycles for Lefschetz type singularities are as follows: Going along the line segment that connects the two singularities, as one approaches the singularity on the left, the cycle $a+b$ vanishes, and as one approaches the singularity on the right, the cycle $c-b$ vanishes.

Now, we will pass to another kind of deformation which is different in nature from the ones that are described above. Note that for a general smooth map $F\colon\f{R}^4\to\f{R}^2$, the differential $dF_p \colon\f{R}^4 \to \f{R}^2$ at a critical point $p$ can have rank either $0$ or $1$. If the rank is $1$, then around $p$ we can find local coordinates such that $F$ is of the form $(t,x,y,z)\to (t, f(t,x,y,z))$ by the inverse function theorem. Similarly, any perturbation $F_s$ of $F$ around $p$ can be expressed in the form $(t,x,y,z)\to(t,f_s(t,x,y,z))$. Therefore, the above moves involved the case where the deformation is focused around a critical point $p$ of $F$ such that $dF$ has rank $1$. In the case of a wrinkled fibration, these are precisely the points lying in the $1$--dimensional part of the critical point set. In fact, any generic deformation around such a critical point is given by one of the above deformations in some coordinate chart. We will elaborate more on this point in the next section using techniques from singularity theory. Our next move will be deforming $F$ around a point $p$ such that $dF_p$ vanishes. For our purposes, these correspond to deforming a wrinkled fibration around a Lefschetz type singularity. 

\paragraph{Move 4 (Wrinkling) :}  Around a Lefschetz type singularity, we have oriented charts where $F\colon\f{R}^4 \to \f{R}^2 $ is given by $(t,x,y,z) \to (t^2-x^2+y^2-z^2, 2tx+2yz)$, or in complex coordinates $u=t+ix$ and $v=y+iz$, $F$ is given by $(u,v)\to u^2+v^2 $. Now the simplest non-trivial deformation of such a map is given by the map $F_s\colon\f{C}^2\to \f{C}$ defined by 
\[ (u,v)\to u^2+v^2+ s \text{Re}u\]

or in real coordinates:
\[ (t,x,y,z) \to (t^2-x^2+y^2-z^2+st, 2tx+2yz) \]
The stability of this map follows from a standard result in singularity theory, see Morin (\cite{morin}). Therefore, the family $F_s$, for $s\in[0,1]$, indeed gives us a family of wrinkled fibrations. The critical points of $F_s$ are the solutions of $x^2+t^2+\frac{st}{2}=0, y=z=0$. This circle can be parametrized by $t=-\frac{s}{4}(1+\cos\theta), x=\frac{s}{4}\sin\theta$, and the critical value set is given by  $ \{(-\frac{s^2}{8}(1+\cos\theta)(2-\cos\theta), -\frac{s^2}{8}(1+\cos\theta)\sin\theta): \theta \in [0,2\pi]\}$. It is easily checked that this equation  defines a curve with 3 cusps. $F_0$ is the standard map around a Lefschetz type singularity, and $F_s$ for $s>0$ is a wrinkled fibration with $3$ wrinkles as shown in Figure $6$. We will refer to the critical value set of this map as triple cuspoid.

\begin{figure}[!ht]
\begin{center}
\includegraphics{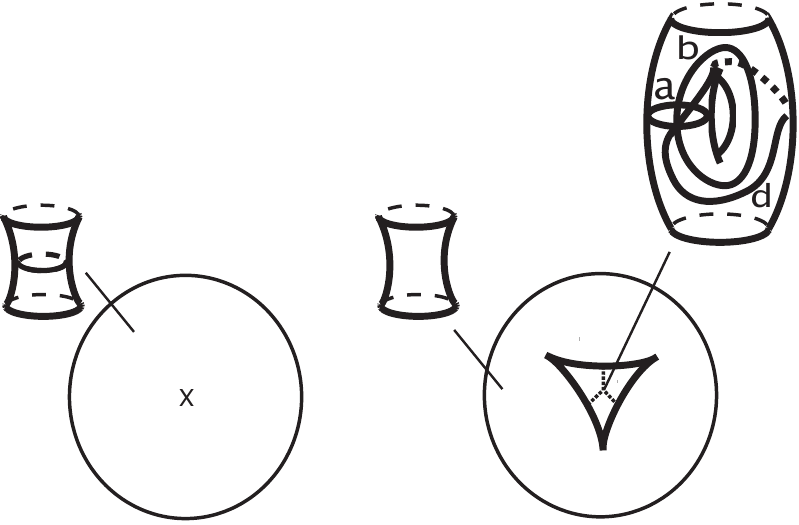}
\end{center}
\caption{Wrinkling}
\end{figure}

The vanishing cycles are $a$, $b$ and $d = b+c$, where $a$,$b$ and $d$ are depicted on the right side of Figure $6$. The curves $a$, $b$ and $c$, which also appear in the middle picture of Figure $5$, are taken to be the standard set of generators for the doubly punctured torus. As shown in Figure $6$, we can in fact arrange so that $d$ passes through the intersection point of $a$ and $b$ and intersects $a$ and $b$ transversely at that point.
 
The importance of this configuration is that all three cycles intersect at a point transversely and there is no path connecting the two boundary components of the doubly punctured torus that does not intersect these three cycles. More precisely, given a configuration of $3$ simple closed curves on a doubly punctured torus with this property, there is a diffeomorphism of the doubly punctured torus which brings the set of curves to the curves $a$, $b$ and $d$ as in Figure $6$ ($d$ is a simple closed curve that is homologous to $b+c$ and passes through the intersection point of $a$ and $b$).

A way to see that the vanishing cycles are as claimed is by considering the fibre above a point $w$ as a double covering of $\f{C}$ branched along $2$ or $4$ points depending on whether $w$ lies outside of the triple cuspoid or in the interior region bounded by the triple cuspoid. Specifically, the fibre above $w$ is given by $v^2=w-u^2-s \text{Re}u$, and projecting to the $u$ component gives a double cover of $\f{C}$ branched along $\{u\in \f{C} \colon u^2+s\text{Re}u=w\}$. 
Let $w=w_1+iw_2$, then in real coordinates one can express the branch locus as: 
\begin{eqnarray*} 
  t^2-x^2+st & = & w_1 \\ 
  2tx & = & w_2
\end{eqnarray*}
For the rest of the argument, assume for simplicity that $s=2$. Take a regular value lying in the interior region of the critical values of $F_s$, such as $w = (w_1,w_2) = ( -1/2 , 0 )$. Connect this to the exterior by the arc of points $(-k/2,0)$, $k \in [1,3]$. One can calculate that the branch points are given by either $t =0$, and $x= \pm \sqrt{k/2}$, or  $x =0$, and $t = -1 \pm \sqrt{4-2k} /2$ . Note that, when $k< 2$ , we have four branch points (fibre is double punctured torus) , and when $k > 2$ we have two (fibre is cylinder). The change is the first two branch points corresponding to $t=0$ more or less stay the same, whereas the branch points corresponding to $x=0$ come together along a segment and disappear when $k >2$.

\begin{figure}[!ht]
\begin{center}
\includegraphics{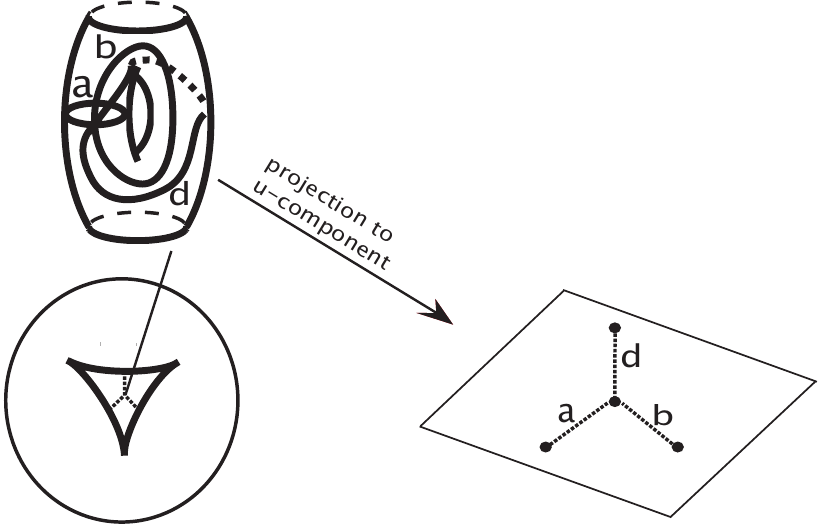}
\end{center}
\caption{The fibre as a double branched cover}
\end{figure}

To get the other vanishing cycles one has to vary $w$ in other directions. The second one can be obtained by $w_1 = -1/2$ and $w_2 = 2k$ where $k$ goes from $0$ to $1$ and the third one can be obtained by $w_1 = -1/2$ and $w_2 = 2k$ where $k$ goes from $0$ to $-1$. One can then see that depending on $k$ we get $4$ branch points if we are in the interior region of the critical values or we have $2$ solutions if we are in the exterior. Corresponding to each of the two variations as above, there are two points in the branch locus which come together whereas the other two stay more or less the same. More precisely, one can verify that corresponding to each direction, the four branch locus points collapse either along $a$, $b$ or $d$ as described in Figure $7$. 

The preimages of these paths by the branched covering map are precisely the vanishing cycles which were also denoted by $a$, $b$ and $d$ on the doubly punctured torus (Figure $7$).  Hence one concludes that the three vanishing cycles intersect transversely at a point. Moreover, it is easy to see by explicit calculation as above that as one approaches a cuspidal point for the fibration $F_s$, in the branched covering picture three of the four branch points come together. For example, if the vanishing cycles $a$ and $b$ collapse as one approaches a cusp singularity of $F_s$, then the end points of the paths $a$ and $b$ come together in the base of the branched cover picture. Reversing our viewpoint, as one crosses a cusp singularity from the low-genus side to the high-genus side the topology of the fibres of $F_s$ is modified by a surgery in a neighborhood of a point in the fibre, which is the preimage of one of the two branch points of the double branched covering map. More precisely, the surgery that we mean is removing a tubular disc neighborhood of a point and gluing back in a punctured torus. We will use this important observation in the next paragraph. 
\paragraph{Deformation of a wrinkled fibration to a broken fibration}\  

Now, we are ready to prove that the local modification of Section \ref{section2} can be obtained by a combination of merging, flipping and wrinkling deformations. Therefore, as promised the local modification given in Section \ref{section2} is also a deformation of wrinkled fibrations. The outcome of this paragraph is the statement that any wrinkled fibration can be deformed to obtain a genuine broken fibration.

\begin{figure}[!ht]
\begin{center}
\includegraphics{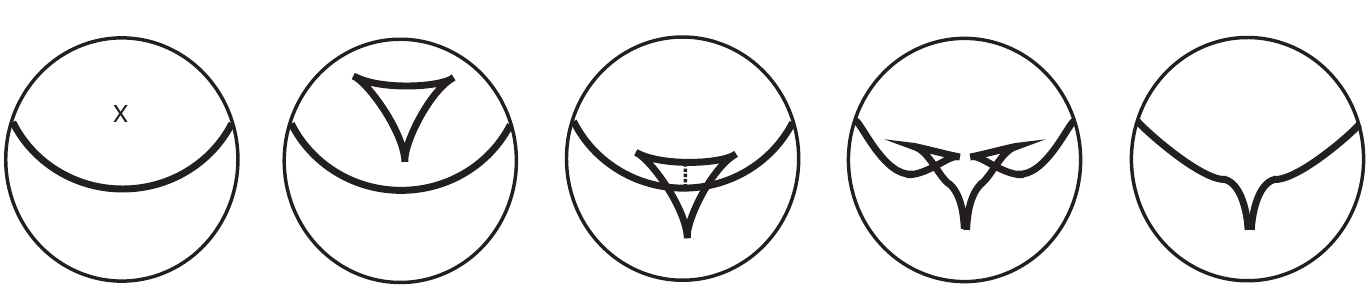}
\end{center}
\caption{Local Deformation}
\end{figure}

Following Figure $8$, first we deform the Lefschetz singularity to a triple wrinkle by applying the wrinkling deformation. Now the key observation here is that we can arrange so that the vanishing cycles corresponding to the bottom cusp of the triple cuspoid do not interfere with the vanishing cycle corresponding to the arc we started with. We will explain this in detail below. Therefore, we can isotope the fibration to the third picture in Figure $8$. Next, we will verify that one can perform a merging move along the dotted line depicted in the third picture in Figure $8$. For this one just needs to verify that the relevant vanishing cycles are in the correct configuration so as to match with the starting point of the local model for the merging move. This will allow us to pass to the fourth picture. Finally, we perform two flipping moves to get to the final result that we wanted.  

Let us now describe the missing pieces of the proof in more detail. First, let's see why one can isotope the second fibration to the third fibration in Figure $8$. For this, we will need to identify various vanishing cycles for the second fibration and observe indeed that the vanishing cycles corresponding to the bottom cusp do not interfere with the vanishing cycle corresponding to the arc. For the fibration that we start with, fix a reference fibre at a point $p$ halfway between the Lefschetz singularity and the singular arc. Recall that the fibre is a punctured torus, and without loss of generality we can assume that the vanishing cycle for the Lefschetz singularity is the $a$ curve and moving towards the arc singularity the $b$ curve vanishes, where $a$ and $b$ are drawn on the left side of in Figure $9$. Now, let's apply the wrinkling move to the Lefschetz singularity. Consider a line segment from $p$ to a central point $q$ of the triple wrinkle passing through a cusp point (drawn as a dotted line on the right side of Figure $9$). As described in the previous section, starting from $p$ if we move along this line segment the fibre above $p$ undergoes a surgery around a neighborhood of a point on the fibre and the genus increases by $1$. Now since the wrinkling move only affects a tubular neighborhood of the curve $a$, after the modification of the Lefschetz singularity by wrinkling move we can choose a reference fibre that is based at the point $q$ which looks like the one drawn in the middle of Figure $9$. In particular, the part of the fibre above $p$ outside of the tubular neighborhood of $a$ is canonically identified to the part of the fibre above $q$ outside the doubly punctured torus that appeared after surgery. More importantly, this latter surgery occurs around a neighborhood of a point which can be isotoped (if necessary) to be disjoint from the $b$ curve. Hence one can parallel transport the $b$ curve from the fibre above $p$ to the fibre above $q$, since the place where the surgery occurs is disjoint form the curve $b$. In particular, the image of $b$ in the fibre above $q$ is disjoint from the vanishing cycles that correspond to the cusp singularity, which are two simple closed curves on the doubly punctured torus which intersect transversely, we denote them by $\alpha$ and $\beta$. Therefore, by applying a diffeomorphism of the doubly punctured torus if necessary the reference fibre above $q$ can be chosen as shown on the right side of Figure $9$. Now, it is clear that one can isotope the second fibration to the third fibration in Figure $8$, since the vanishing cycle $b$ is disjoint from $\alpha$ and $\beta$.

\begin{figure}[!ht]
\begin{center}
\includegraphics{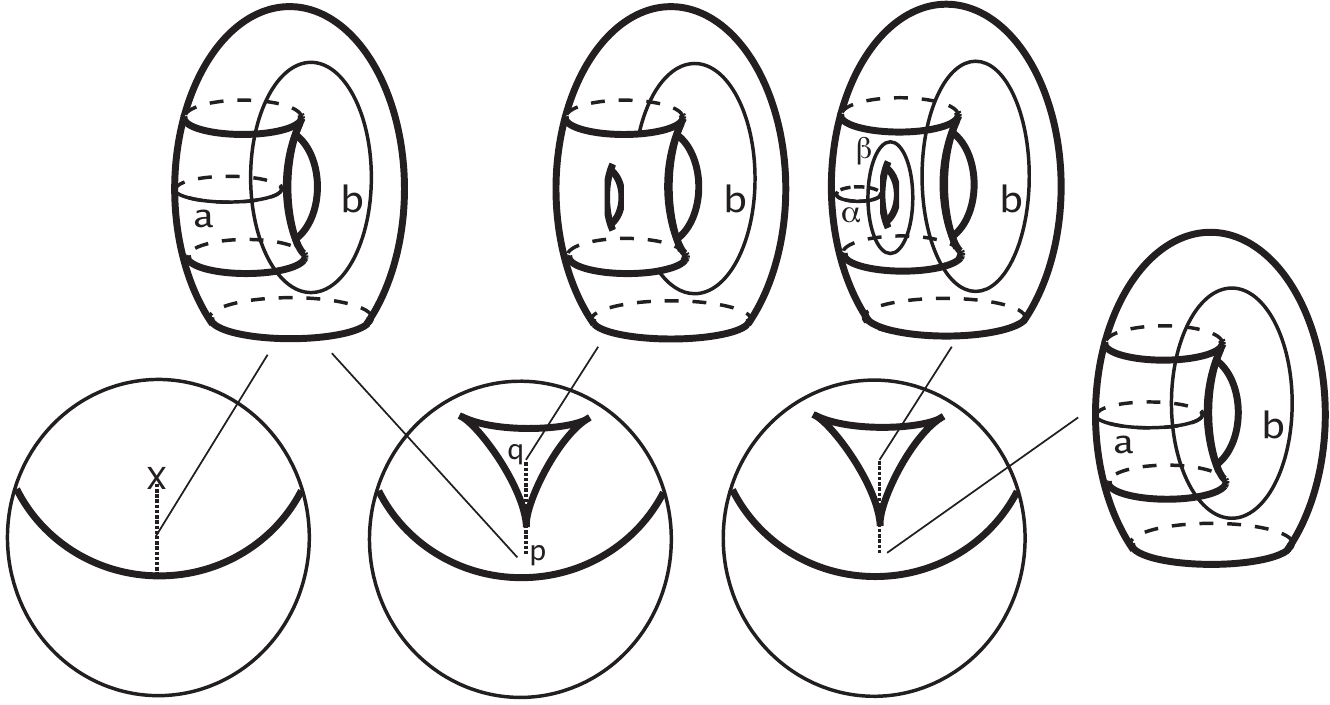}
\end{center}
\caption{Reference fibres}
\end{figure}

Next, to pass from the third fibration to the fourth fibration in Figure $8$, we use a merging move. In order to do that, we need to understand the vanishing cycles above the dotted line segment in the third picture in Figure $8$. Choose a reference fibre above a point in the middle of the dotted line segment. As before, we can standardize it so that it looks like the right side of Figure $9$. Now, as one goes down the curve $b$ vanishes and as one goes up the vanishing cycle $\gamma$ has the properties that it lies in the doubly punctured torus,  intersects $\alpha$ and $\beta$ at their intersection point and any path connecting the boundary circles of the doubly punctured torus has to intersect the union of $\alpha$, $\beta$ and $\gamma$. Therefore, comparing Figure $9$ with Figure $6$, $b$ has to intersect $\gamma$ once. Hence we can perform a merging move. 

Finally, we apply two flipping moves to the fourth fibration in Figure $8$ to pass to the fifth fibration. These are also allowed, since the configuration of $\alpha$, $\gamma$, $b$ and the configuration of $\beta$, $\gamma$, $b$ match the configuration of $a$, $b$, $c$ in Figure $5$ of the flipping move. This completes the proof of the fact that the fibration on the left of Figure $8$ is a deformation of the fibration on the right.

\section{Generic deformations of wrinkled fibrations and $(1,1)$--stability}
\label{section4}
In this section, we prove that the set of moves listed in the Section \ref{section3} are sufficient to produce any deformation of wrinkled fibrations. More precisely, we prove the following theorem:
\begin{theorem}
\label{uniqueness}	
	Let $X$ be a compact $4$--manifold, and let $F_s\colon X \to \Sigma$ be a deformation of wrinkled fibrations. Then it is possible to deform $F_0$ to $F_1$ by applying to $F_0$ a sequence of the four moves described in Section \ref{section3} and isotopies staying within the class of wrinkled fibrations.
			
\end{theorem}

\emph{Proof.}
First, observe that we can get rid of the Lefschetz type singularities of $F_0$ and $F_1$ using the wrinkling move. So we can assume that $F_0$ and $F_1$ have no Lefschetz type singularities. Also, since Lefschetz type singularities are unstable under small deformations, we can assume that the deformation does not create any new Lefschetz singularity. More precisely, we perturb the deformation by keeping the end points fixed so that we avoid any creation of critical points where $dF_s$ vanishes. This is possible since purely wrinkled fibrations are stable under small perturbation whereas the existence of points where $dF_s$ vanishes is not generic. Therefore, we have reduced to the case where $F_s$ is a wrinkled fibration except for finitely many values of $s$ such that $F_s$ has no Lefschetz singularity for all $s$. So, we can assume that around a critical point $p$ of $F_{s_0}$ for any $s_0 \in [0,1]$, we have coordinate charts so that for $s\in [s_0-\epsilon,s_0+\epsilon]$, $F_s \colon \f{R}^4 \to \f{R}^2 $ is given by $(t,x,y,z)\to (t,f_s(t,x,y,z))$ and $f_{s_0}(0)=df_{s_0}(0)=0$. We will next show that generically $f_s$ is given by one of the $3$ models described in Section \ref{section3} corresponding to the moves birth, merging and flipping. For this, we will introduce the notion of \emph{$(1,1)$--stable unfoldings} following Wasserman \cite{wasserman} and give a classification of such maps using the machinery developed in \cite{wasserman}, which in turn is based on the celebrated classification of unfoldings by Thom.\\ 

\begin{definition}
\label{defrsequiv}
Let $f\colon\f{R}^5\to \f{R}$ and $g\colon\f{R}^5\to \f{R}$ be map germs with $f(0)=g(0)=0$. With $f$ we associate a germ $F\colon\f{R}^5\to\f{R}^3$, defined by $F(s,t,x,y,z)=(s,t,f(s,t,x,y,z))$. Similarly we associate a germ $G\colon\f{R}^5\to\f{R}^3$ with $g$, given by $G(s,t,x,y,z)=(s,t,g(s,t,x,y,z))$. We say that $f$ and $g$ are $(1,1)$--equivalent if there are germs at $0$, $\Phi \in \text{Diff}(\f{R}^5), \Lambda \in \text{Diff}(\f{R}^3)$ and $\psi\in \text{Diff}(\f{R}^2)$, and $\phi \in \text{Diff}(\f{R})$ fixing the origin such that the following diagram commutes:
\begin{displaymath}
\xymatrix{
\f{R}^5 \ar[r]^F \ar[d]_\Phi & \f{R}^3 \ar[r]^p \ar[d]_\Lambda & \f{R}^2 \ar[r]^q \ar[d]_\psi & \f{R} \ar[d]_\phi \\
\f{R}^5 \ar[r]^G & \f{R}^3 \ar[r]^p & \f{R}^2 \ar[r]^q & \f{R} }
\end{displaymath}

where $p\colon\f{R}^3\to \f{R}^2$ is the projection onto the first factor and $q\colon\f{R}^2\to\f{R}$ is the projection onto the first factor. 	
\end{definition}

Note that if a one-parameter family deformation $F_s$ of wrinkled fibrations is represented by $(t,x,y,z)\to$ $(t, f(s,t,x,y,z))$ in some coordinate charts and $g$ is $(1,1)$--equivalent to $f$, then we can find coordinate charts such that the deformation is represented by $(t,x,y,z)\to(t,g(s,t,x,y,z))$ in these new coordinate charts. Therefore, in order to complete the proof of theorem \ref{uniqueness}, we need a classification theorem of \emph{generic} functions up to $(1,1)$--equivalence, which we state after making precise what generic means.

\begin{definition}
\label{defstable}
Let $\mathcal{E}(\f{R}^n,\f{R}^p)=$ set of germs at $0$ of smooth mappings from $\f{R}^n\to \f{R}^p$. Let $\mathcal{E}(s,t,x,y,z)=\mathcal{E}(\f{R}^5,\f{R}), \mathcal{E}(s,t)=\mathcal{E}(\f{R}^2,\f{R}), \mathcal{E}(s)=\mathcal{E}(\f{R},\f{R})$ such that the labels reflect the parameters that we are using.

Let $f\colon\f{R}^5\to \f{R}$ with $f(0)=0$ and let $F\colon\f{R}^5\to\f{R}^3$ be given by $F(s,t,x,y,z)=(s,t,f(s,t,x,y,z))$. We say that $f$ is infinitesimally  $(1,1)$--stable if \[ \mathcal{E}(s,t,x,y,z)=\langle \frac{\del f}{\del x} , \frac{\del f}{\del y}, \frac{\del f}{\del z} \rangle {\mathcal{E}(s,t,x,y,z)} + \langle \frac{\del f}{\del t} \rangle {\mathcal{E}(s,t)} + \langle \frac{\del f}{\del s} \rangle {\mathcal{E}(s)} + F^*{\mathcal{E}(\f{R}^3,\f{R})}. \]

\end{definition}

One may interpret this condition geometrically as saying roughly that the ``tangent space'' at $f$ to the $(1,1)$--equivalence class of $f$ is maximal, i.e. is equal to the ``tangent space'' to the unique maximal ideal in $\mathcal{E}(\f{R}^5, \f{R})$ consisting of the set of germs $f$ such that $f(0)=0$.

We remark here that by Theorem $3.15$ in \cite{wasserman} any perturbation of a $(1,1)$--stable germ in weak $C^\infty$--topology can be represented by a $(1,1)$--stable germ. Therefore, in this sense, a generic deformation will be $(1,1)$--stable. \\

\begin{theorem} 
\label{classification}
Let $f\colon\f{R}^5\to\f{R}$ be a $(1,1)$--stable germ with $f(0)=0$. Then $f$ is $(1,1)$--equivalent to one of the following functions as germs:
\begin{eqnarray*} 
 h_0(s,t,x,y,z)& = & \pm x^2 \pm y^2 \pm z^2 \ \ \ \ \ \ \ \ \ \ \ \ \ \ \text{Morse singularity} \\ 
 h_1(s,t,x,y,z) & = & x^3+tx \pm y^2 \pm z^2 \ \ \ \ \ \ \ \ \ \ \text{Cusp singularity} \\ 
h_2(s,t,x,y,z) & = & x^3+t^2x+sx \pm y^2 \pm z^2 \ \ \text{Birth move} \\
h_3(s,t,x,y,z) & = & x^3-t^2x+sx \pm y^2 \pm z^2  \ \ \text{Merging move} \\
h_4(s,t,x,y,z) & = &  x^4+x^2s+xt \pm y^2 \pm z^2  \ \ \text{Flipping move}
\end{eqnarray*} 

\end{theorem}

The proof of Theorem \ref{classification} is given in the Appendix A. This completes the proof of Theorem \ref{uniqueness} where the signs in the statement of Theorem \ref{classification} are determined by imposing the condition that the maps become wrinkled fibrations. 

\section{The corresponding deformations on near-symplectic manifolds}
\label{section5}
\begin{theorem}
\label{canonicalsymp}
Let $X$ be a compact $4$--manifold, and let $f\colon X\backslash P \to \Sigma$ be a wrinkled pencil. Let $Z$ denote the $1$--dimensional part of the critical value set of $f$. Suppose that there exists a cohomology class $h \in H^2(X)$ such that $h(F)>0$ for every component $F$ of every fibre of $f$, then there exist a near-symplectic form $\omega$ on $X$, with zero set $Z$ and such that $\omega$ restricts to a symplectic form on the smooth fibres of the fibration. Moreover, $\omega$ determines a deformation class of near-symplectic forms canonically associated to $f$.
\end{theorem}

Note that if every component of every fibre of $f$ contains a point in $P$, then the cohomological assumption holds automatically. We will not give a full proof of this theorem as the proof in \cite{adk} applies here almost verbatim. The only modification required is in part $1$ of the proof given in \cite{adk}, where one constructs a near-symplectic form positive on the fibres which is defined only in a neighborhood of the critical point set. For the wrinkled fibrations, we introduce a new type of singularity on the critical value set, namely the cusp singularity. Therefore one needs to say a word about how to construct a near-symplectic form positive on the fibres for the local model of the cusp singularity. For that, recall the local model for the cusp singularity. To wit, we have oriented charts where the wrinkled fibration is given by: \[ f\colon (t,x,y,z) \to (t,x^3-3xt+y^2-z^2) \] Now, consider the $2$--form $\omega = dt \wedge df_t + * (dt \wedge df_t) $, where $f_t(x,y,z)=x^3-3xt+y^2-z^2$ are Morse except at $t=0$. This form is self-dual by construction. Since $f_t$ is Morse except at $t=0$, this form is transverse to the $0$--section of $\Lambda^+$. The only missing property for $\omega $ to be near-symplectic is that it be closed. In fact, in this specific example of $f_t$ that we are considering $\omega $ is not closed. The reason that we are considering this specific $\omega$ is because it is positive on the fibres by construction. Therefore, we want to modify $\omega$ by adding some terms so that it is closed and at the same time preserve the property that it is positive on the fibres. In this section, this will be the general scheme for finding explicit near-symplectic forms on a given fibration. One such modification is as follows:
\[ \tilde{\omega}= dt \wedge df_t + * (dt \wedge df_t)- y(3dt\wedge dz + 6x dz\wedge dx)\] However, in order to control the positivity we need to ensure that the extra terms we added are small when evaluated on a basis of a fibre. Therefore, we multiply that term with an $\epsilon>0$, and in order to have a closed form we need to also multiply the $dx \wedge dt+dy \wedge dz$ component of $dt \wedge df_t + * (dt \wedge df_t)$ also by $\epsilon$. In what follows, we will do this modification several times, therefore we introduce a scaling map $R_\epsilon \colon \Omega^2_+ \to \Omega^2_+ $ given by:
\begin{eqnarray*}
R_\epsilon(dt\wedge dx + dy \wedge dz) & = & \epsilon(dt\wedge dx + dy \wedge dz) \\
R_\epsilon(dt\wedge dy + dz \wedge dx) & = & (dt\wedge dy + dz \wedge dx) \\
R_\epsilon(dt\wedge dz + dx \wedge dy) & = & (dt\wedge dz + dx \wedge dy) 	
\end{eqnarray*}
 So finally we have our near-symplectic form given by: 
\begin{eqnarray*}
\omega_\epsilon & = & R_\epsilon(dt \wedge df_t + * (dt \wedge df_t))- \epsilon y(3dt\wedge dz + 6x dz\wedge dx) \\
& = &  3\epsilon (x^2-t) (dt\wedge dx +dy\wedge dz) \\
& + & 2y dt\wedge dy+(2y-6\epsilon x y ) dz\wedge dx \\
& - & (2z+3\epsilon y) dt \wedge dz -2z dx\wedge dy 
\end{eqnarray*}

Now, choose $\epsilon \leq 1/6$. Then one can check easily that $\omega_\epsilon$ is a near-symplectic form on $D^4$ and its restriction to smooth fibres of $f$ are symplectic. Thus, we can use $\omega_\epsilon$ for the local construction in the proof of Theorem \ref{canonicalsymp}.

Theorem \ref{canonicalsymp} tells us that there is a natural deformation class of near-symplectic forms on each of the local models of wrinkled fibrations. In what follows, we will give explicit models of near-symplectic forms for each of the local model of wrinkled fibrations described in the previous sections. Furthermore, we will provide one-parameter families for the deformations corresponding to the $4$ moves given in Section \ref{section3}. These will be near-symplectic cobordisms in the sense of the following definition given by Perutz \cite{perutz}.

\begin{definition} 
	A one-parameter family $\{\omega_s\}_{s\in[a,b]}$ of closed $2$--forms on $X$ is called a near-symplectic cobordism if, for all $(x,s)\in X\times [a,b]$, either $(\omega_s \wedge \omega_s) (x) > 0 $ or, $\omega _s(x)=0$ and $(\nabla\omega)(x,s)$ has rank $3$. 
	
\end{definition}

The strategy will be the same as the construction of the local model around a cusp singularity. We first exhibit a $2$--form positive on the fibres which is not necessarily closed. Then we modify it by adding small terms. We will mostly restrict the domain of the wrinkled fibration to $D^4$ to ensure positivity. Since every deformation is local and the critical value set lies in $D^4$, this is not different from the previous considerations. 

\paragraph{Deformation 1 (Birth) :} The deformation is given by 
$F_s \colon D^4 \to \f{R}^2$ :
\[ (t,x,y,z) \to (t, x^3+3(t^2-s)x+y^2-z^2)\]
Let $f_s = x^3+3(t^2-s)+y^2-z^2$. Consider the deformation: \begin{equation} \omega_s = R_\epsilon(dt \wedge df_s + *(dt \wedge df_s)) + 6 \epsilon y (t dt \wedge dz + x dx \wedge dz )  \end{equation}

This form is closed and if we choose $\epsilon\leq 1/6$, it is near-symplectic on $D^4$. Furthermore, an easy calculation shows that $\omega_s$ is symplectic on smooth fibres of $F_s$. Now, here we remark that $\omega_s$ is in fact precisely the Luttinger--Simpson model of birth of a circle singularity which was defined in the introduction to this paper. Therefore, the maps $F_s$ gives a family of wrinkled fibrations adapted to the model of Luttinger--Simpson of near-symplectic cobordism $\omega_s$.

\paragraph{Deformation 2 (Merging) :} The deformation is given by 
$F_s \colon D^4 \to \f{R}^2$ :
\[ (t,x,y,z) \to (t, x^3+3(s-t^2)x+y^2-z^2)\]
Let $f_s = x^3+3(s-t^2)+y^2-z^2$. Consider the deformation: \begin{equation} \omega_s = R_\epsilon(dt \wedge df_s + *(dt \wedge df_s)) - 6 \epsilon y (t dt \wedge dz + x dz \wedge dx ) \end{equation}

As before, this form is closed and for $\epsilon \leq 1/6$, it is near-symplectic on $D^4$. This is a variation of the birth model, the zero-set undergoes a surgery by addition of a one handle. Again, this is a near-symplectic cobordism, and the family $F_s$ is adapted to $\omega_s$, i.e., the restriction of $\omega_s$ to smooth fibres of $F_s$ is positive.

\paragraph{Deformation 3 (Flipping) :} We again follow the same strategy as above. However, in this case we do not need to restrict to $D^4$. Namely, consider the deformation for flipping move given by $F_s \colon \f{R}\times \f{R}^3 \to \f{R}^2$ :
\[ (t,x,y,z) \to (t, x^4-x^2 s+xt+y^2-z^2)\]

Now, let $f_s = x^4-x^2 s +xt+y^2-z^2$. Then we calculate: 
\begin{eqnarray*} 
dt\wedge df_s + *(dt\wedge df_s) & = & (4x^3-2xs+t)(dt \wedge dx +dy \wedge dz) \\
 & + & 2y (dt \wedge dy + dz \wedge dx)  \\ 
 & - & 2z (dt \wedge dz + dx \wedge dy) 
\end{eqnarray*}
This form is positive when restricted to the smooth fibres of $F_s$ by design. However, this form is not closed. Therefore, to make it closed we modify it naively as follows to make it closed. 
\begin{eqnarray} 
\nonumber \omega_s & = & (4x^3-2xs+t)( dt \wedge dx + dy \wedge dz) \\
 & + & (2y-2z) dt \wedge dy + (12x^2-2s+2)y dz \wedge dx \\ 
\nonumber & - & (2z+y) dt \wedge dz - (12x^2-2s+1)2z dx \wedge dy 
\end{eqnarray}
Now $\omega_s$ is closed and in fact an easy calculation shows that for $s \leq 1/3$, $\omega_s$ is still positive when restricted to the smooth fibres of $F_s$. Furthermore, the zero locus of $\omega_s$ is exactly the critical point set of $F_s$. Therefore, we conclude that $\omega_s$ in fact belongs to the canonical class of near-symplectic forms provided by Theorem \ref{canonicalsymp} for the fibration $F_s$. Furthermore, the near-symplectic cobordism $\omega_s$ for $s\in[-1,1/3]$ is through near-symplectic forms, that is, for each $s\in[-1,1/3]$, $\omega_s$ is near-symplectic and adapted to $F_s$ in the sense of Theorem \ref{canonicalsymp}. Hence, we conclude that the flipping move does not alter the near-symplectic geometry.

\paragraph{Deformation 4 (Wrinkling) :}Recall that the wrinkling move is given by $F_s\colon\f{R}^4 \to \f{R}^2$ : 
\[(t,x,y,z) \to (t^2+st-x^2+y^2-z^2, 2tx+2yz)\] Let $f_s = t^2+st-x^2+y^2-z^2$ and $g=2tx+2yz$. Then a natural candidate for an adapted near-symplectic form for $F_s$ is given by $df_s\wedge dg + *(df_s\wedge dg)$ However, as before, 
\begin{eqnarray*}
*(df_s \wedge dg) & = &((2t+s)2t+4x^2) dy \wedge dz + (4y^2+4z^2) dt\wedge dx \\
& + & ((2t+s)2z-4xy) dz \wedge dx + (4xy - 4tz) dt \wedge dy \\
& + & ((2t+s)2y+4xz) dx \wedge dy - (4xz +4ty) dt \wedge dz 
\end{eqnarray*} 

is not closed. Therefore we modify it to the following form. 
\begin{eqnarray*}
\sigma_s & = &((2t+s)2t+4x^2) dy \wedge dz + (4y^2+4z^2) dt\wedge dx \\
& + & 2 ((2t+s)2z-4xy) dz \wedge dx + 2(4xy - 4tz -sz)dt \wedge dy 
\end{eqnarray*}
It is an easy calculation to check that $\sigma_s$ is closed and positive when restricted to the fibres of $F_s$. Now, in order to get a near-symplectic form, we restrict to $D^4$, so $F_s \colon D^4 \to D^2$, and to $\sigma_s$ we add a large multiple of the pullback of the standard symplectic form on $D^2$ by $F_s$. Thus,
\begin{equation}  \omega_s = k (df_s \wedge dg) + \sigma_s 
\end{equation}
for $k$ large enough is an adapted near-symplectic form, that is, it vanishes exactly at the critical value set of $F_s$ and restricts positively to smooth fibres of $F_s$. Observe that, here also we can see a birth of a zero-circle happens as $s$ goes through negative values to positive values. Therefore, it is possible that this form is deformation equivalent through near-symplectic forms to Luttinger--Simpson model.
\section{Applications}
\label{section6}
\textbf{Merging of zero-sets:}
Here we reprove the Theorem 1.4 in \cite{perutz} using moves on broken fibrations. 

\begin{theorem}
\label{6.1}

Given a connected near-symplectic manifold $(X,\omega_0)$, with $\omega_0$ having a zero-set with $n$ components, where $n \geq 1$, one can find a near-symplectic cobordism $\omega_{[0,1]}$ such that $\omega_1$ has $k$ components for any given $k \geq 1$. Furthermore, this near-symplectic cobordism is equipped with an adapted wrinkled pencil.

\end{theorem}

Our proof will be obtained by applying moves on a broken pencil adapted to the given near-symplectic manifold. However, one can ignore the base points of the pencil since all the modifications will take place away from them. In this way, we obtain a quicker proof as well as our deformation includes a deformation of wrinkled fibrations associated to it. 

\emph{Proof.} Choose an adapted broken pencil for $(X,\omega_0)$ which exists by the main construction in \cite{adk}. The proof is divided into two parts according to increasing or decreasing the number of components of the zero-set. 

First, let's show that we can add a new component. Restrict the given pencil to a smooth $D^2$ fibration over $D^2$, which is isolated from the singularities of the broken pencil and apply the birth move. Deformation $1$ above, tells us that this gives us a near-symplectic cobordism $\omega_{[0,1]}$, where $\omega_1$ has one more component in its zero-set. 

Second, let's show that if $n > 1$, we can find a near-symplectic cobordism where the number of components decreases by $1$. This part will be longer, since we can't directly apply the merging move as the configuration needed for the merging move is not always possible to achieve. However, we will apply an alternative combination of moves to produce a merging of zero-components in the total space. First, choose two distinct components of the zero-set. Now, connect these components by a path $\alpha\colon[0,1] \to X$ such that the following properties are satisfied. 
\begin{itemize}
\item $\alpha^{-1}(Z)=\{0,1\}$ where $Z$ is the zero set of $\omega_0$.
\item $\alpha'(0), \alpha'(1) \in L^+$, where $L^+$ is the positive eigen-subbundle of $NZ$ as defined in the introduction. 
\item $\alpha$ is transversal to the fibres of the broken pencil.
\end{itemize}

Clearly, such paths are in abundance. Indeed, locally near the end points it's easy to build the path using the local models; and everywhere else, being transverse to the fibres is generic. Restrict the pencil to a neighborhood $N=U \cup V \cup W$ of $\alpha$, where $U$ and $W$ are preimages of a small neighborhood of the image of $\alpha(0)$ and $\alpha(1)$, and $V$ is a tubular neighborhood $\alpha$. Then we have a picture as depicted on the left of Figure $10$, where the fibres depicted lie in $U$ and $W$. The preimage of the middle region is $V$, at each fibre this cuts out a disc. Now, we can apply two flipping moves to both sides, and obtain a fibration as depicted in the middle part of Figure $10$. 

\begin{figure}[!h]
\begin{center}
\includegraphics{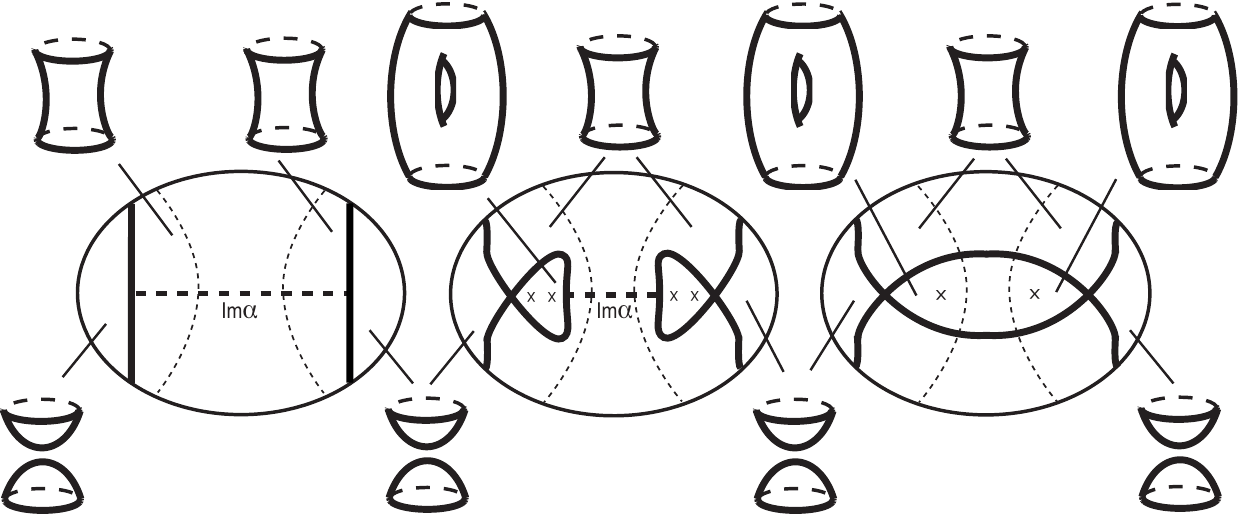}
\end{center}
\caption{Merging of zero-circles along the path $\alpha$}	 
\end{figure}

Notice that these flipping moves do not alter the deformation class of the near-symplectic form and hence the isotopy class of the zero-set is unchanged after these moves, only the broken pencil structure has been changed. Finally, given such a configuration, we can apply an inverse merging move to the fibration (See the remark at the end of the description of the merging move in Section \ref{section3}). In the total space this corresponds to merging of the zero-sets and the deformation of the near-symplectic form is given in the form of a near-symplectic cobordism as in Deformation 2, given by the formula $2$, except $s$ must be replaced by $-s$, as we apply an inverse merging move. 
\paragraph{Broken fibrations with connected fibres:}
Another application of the techniques discussed in this paper is based on an idea of Baykur and also appears in \cite{baykur}. Here we reconstruct that argument for the sake of completeness.

\begin{theorem} 
\label{connectedness}	
	Given a connected near-symplectic manifold $(X,\omega)$, one can always find a broken pencil $f:X\to S^2$, adapted to $\omega$, the fibres of which are connected. 	
	
\end{theorem}

Therefore, in order to define Perutz's Lagrangian matching invariant one can always start with a broken fibration with connected fibres. This indeed simplifies some of discussions in \cite{matching} and allows us to define Lagrangian matching invariant for a slightly larger number of $\text{Spin}^c$ structures.

\emph{Proof:} For simplicity, we start with the case where the zero-set of $\omega$ consists of a single component. Now, observe that, by perturbing $\omega$ away from its zero-set, and using the main result in $\cite{adk}$, we can ensure that there exists an adapted broken pencil for the perturbed near-symplectic form. Since the perturbation can be taken to be arbitrarily small, the latter broken pencil will be adapted to $\omega$ as well. Without loss of generality we can assume that there are no base points, otherwise we blow-up first, apply the argument below and blow-down in the end.

Either the fibres are connected or suppose the fibres above the northern hemisphere have genera $g$, and the fibres above the southern hemisphere have genera $g_1$, $g_2$ such that $g_1+g_2=g$. Since $X$ is connected, this is the only possibility as the fibres above the ``high-genus side'' have to be connected. Furthermore, we can assume that there are no Lefschetz-type singularities in the ``low-genus side'' since they can be isotoped to the ``high-genus side''. This is simply because starting from a regular fibre in the southern hemisphere adjoining a Lefschetz singularity means adding a $2$--handle, and adjoining a broken singularity means adding a $1$--handle. But the order of adding these handles can be reversed by an isotopy, therefore one can first add the $1$--handle corresponding to the broken singularity, then add the $2$--handles corresponding to the Lefschetz-type singularities. Therefore, we can assume that the fibration is trivial above the southern hemisphere. So, the preimage of a neighborhood of the southern hemisphere is given as in Figure $11$.

\begin{figure}[!h]
\begin{center}
\includegraphics{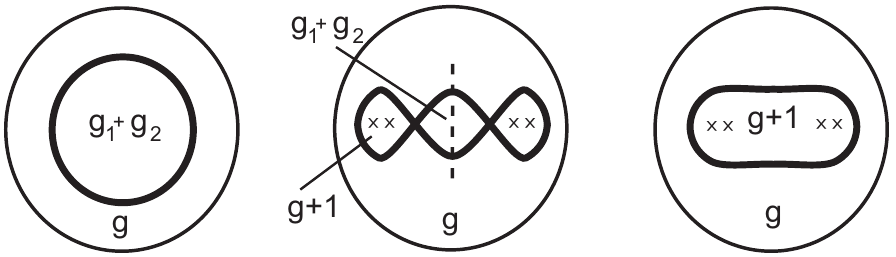}
\end{center}
\caption{Making the fibres connected}	 
\end{figure}

Now, we apply two flipping moves to pass to the middle picture in Figure $11$. Finally, to obtain the final fibration depicted on the right, we perform an isotopy interchanging the two ``legs'' of the flips in the middle. This is allowed, since if we consider an arc cutting these ``legs'' transversely as shown in the middle picture in Figure $11$ as a dotted line, the topology of the fibre changes by first vanishing of a separating cycle (that is, a $2$--handle attachment) and then attaching a $1$--handle. Again, the order of attachment does not matter, hence by an isotopy one can obtain a broken fibration where the fibres are connected.

When the zero locus of $\omega$ consists of more than one circle, these various circles live in disjoint parts of the fibers above the equator of $S^2$. We can again push the Lefschetz fibers to the high genus side (northern hemisphere) and ensure that the fibration is trivial above the southern hemisphere. Since the modification explained above is local in the fibre (it only affects a neighborhood of the vanishing cycle for the broken singularity), it can be performed simultaneously on each of the circles. Pictorially this again amounts to the transition shown in Figure $11$, but with several circles ``stacked'' on top of each other in disjoint parts of the fiber.
\paragraph{Removing achiral Lefschetz singularities:} In this paragraph, we prove that any achiral Lefschetz singularity can be replaced with a circle of broken singularities and three Lefschetz singularities. Recall that an achiral Lefschetz singularity is modeled in orientation preserving charts by the complex map $(w,z)\to \bar{w}^2+z^2$. 

Now, given an achiral Lefschetz singularity we can consider the same deformation that was used in the wrinkling move in Section \ref{section3}. Namely, let $F_s:\f{C}^2\to\f{C}$ given by: \[(w,z)\to \bar{w}^2+z^2+s\text{Re}w\] The map $F_s$ is identical to that considered in Section \ref{section3} up to the orientation-reversing diffeomorphism $(w,z) \to (\bar{w},z)$. Thus its critical values and vanishing cycles are the same as in Section \ref{section3} up to a reversal of the orientation of the fibre. Namely, when $s>0$, we observe a birth of a circle of singularities and we get a wrinkled fibration with $3$ cusps as on the right side of Figure $6$ (except the configuration of the vanishing cycles is reversed). Next, apply the local modification discussed in Section \ref{section2} to replace each of the three cusps in the base by a smooth arc and a Lefschetz singularity. Thus, we have replaced a neighborhood of an achiral singularity with a genuine broken fibration with a new circle of broken singularities together with three Lefschetz singularities.

We remark that the new singular circle obtained here is an even circle, whereas the new singular circle obtained in the original wrinkling move is an odd circle. (The notions of even and odd circles were defined in the introduction.) The fact that the original wrinkling move yields an odd circle follows from the fact that on a near-symplectic manifold the number of even circles is equal to $1-b_1 + b_2^+ $ modulo $2$, as was mentioned in the introduction. A more direct way to see this is as follows: After modifying the cusps as in Section \ref{section2}, we obtain a singular circle of broken singularities and three Lefschetz singularities. Take a small disc including the three Lefschetz singularities but not intersecting the singular circle. Fix a reference fibre above a point on the boundary of this disc such that the curve $a$ vanishes as one approaches the broken singularity. In Figure $7$, after modifying the cusps, the reference fibre we are fixing is on the lower left side of the picture. The monodromy around the boundary of this disc is given by the composition of three right handed Dehn twists corresponding to the Lefschetz singularities. Again, from Figure $7$ and the calculation of vanishing cycles in Section \ref{section2}, one can conclude that this monodromy is given by $\mu=\tau_{a+d}\circ\tau_{b-d}\circ\tau_{a-b}$. From this, we see that $\mu(a)=-a$, which shows that the circle is an odd circle. Now, in the above perturbation for the achiral Lefschetz singularity case, all the configuration is the same except, the base picture in Figure $7$ is reflected so that the counter-clockwise ordering of $a,b,d$ is changed to $a,d,b$. Then, one calculates the monodromy to be \[\mu=\tau_{a+b}\circ\tau_{b+d}\circ\tau_{a-d} \] which gives $\mu(a)=a$. Hence the singular circle obtained is an even circle. 

\begin{corollary} Let $X$ be an arbitrary closed $4$--manifold and let $F$ be a closed surface in $X$ with $F\cdot F = 0$. Then there exists a broken Lefschetz fibration from $X$ to $S^2$ with embedded singular locus, and having $F$ as a fibre. Furthermore, one can arrange so that the singular set on the base consists of circles parallel to the equator with the genera of the fibres in increasing order from one pole to the other.
	
\end{corollary}

\emph{Proof.} The existence follows from Gay and Kirby's theorem \cite{gkirby}, and the above modification of achiral Lefschetz singularities. Let's prove that this can be done in a certain way so that the singular set on the base consists of circles parallel to the equator. First, note that Gay and Kirby's proof places the round singularities on the tropics of Cancer and Capricorn and the ``highest-genus region'' is the annular region between the tropics. Now move any Lefschetz or achiral Lefschetz singularities in the southern hemisphere towards the equatorial region (moving across circles towards the higher genus region as in Theorem \ref{connectedness}). Then in the southern hemisphere we are left with only a bunch of parallel circles on the tropic of Capricorn; the corresponding round $1$-handles all get attached along disjoint braids (i.e. the circle attachments can be stacked on top of each other or commuted). This is what we need to be able to apply the move on Figure $11$ (i.e., on all the circles simultaneously, placing them on top of each other and in different parts of the fibers: first two flipping moves, then an isotopy exactly as in the argument in Theorem \ref{connectedness}. Consequently, we still have circles on the tropic of Capricorn, but oriented in the opposite way (genus increases towards south pole), and some Lefschetz fibers near the south pole (created by the flips in Figure $11$) and the previously given Lefschetz and anti-Lefschetz fibers near the equator. The latter can now be moved towards the north pole by crossing the circles at tropic of Capricorn. Therefore, we can assume that the singular circles are equatorial with all the circles oriented the same way and the  ``highest-genus region'' is over the south pole. Now, we transform one of the achiral singularities to a circle singularity and three Lefschetz singularities. Next, push all the remaining Lefschetz singularities and the achiral singularities (left between the previous circles and the new circle) across the new circle (into the even higher genus region), so the circle is now in equatorial position (parallel to the previous circles). Finally, we repeat this process until there are no more achiral Lefschetz singularities left. 

\newpage
\section{A summary of moves and further questions }
\label{section7}
The table below summarizes our set of moves. Only the base parts of the fibrations are drawn. Each move is drawn in pairs, as a move on wrinkled fibrations and as a move on the corresponding broken fibrations obtained by replacing each wrinkle by an arc together with a Lefschetz type singularity as was discussed in Section \ref{section2}. Also the references to the formulas concerning the changes in the naturally associated near-symplectic forms provided by Theorem \ref{canonicalsymp} are given. 

\begin{figure}[!h]
\begin{center}
\includegraphics{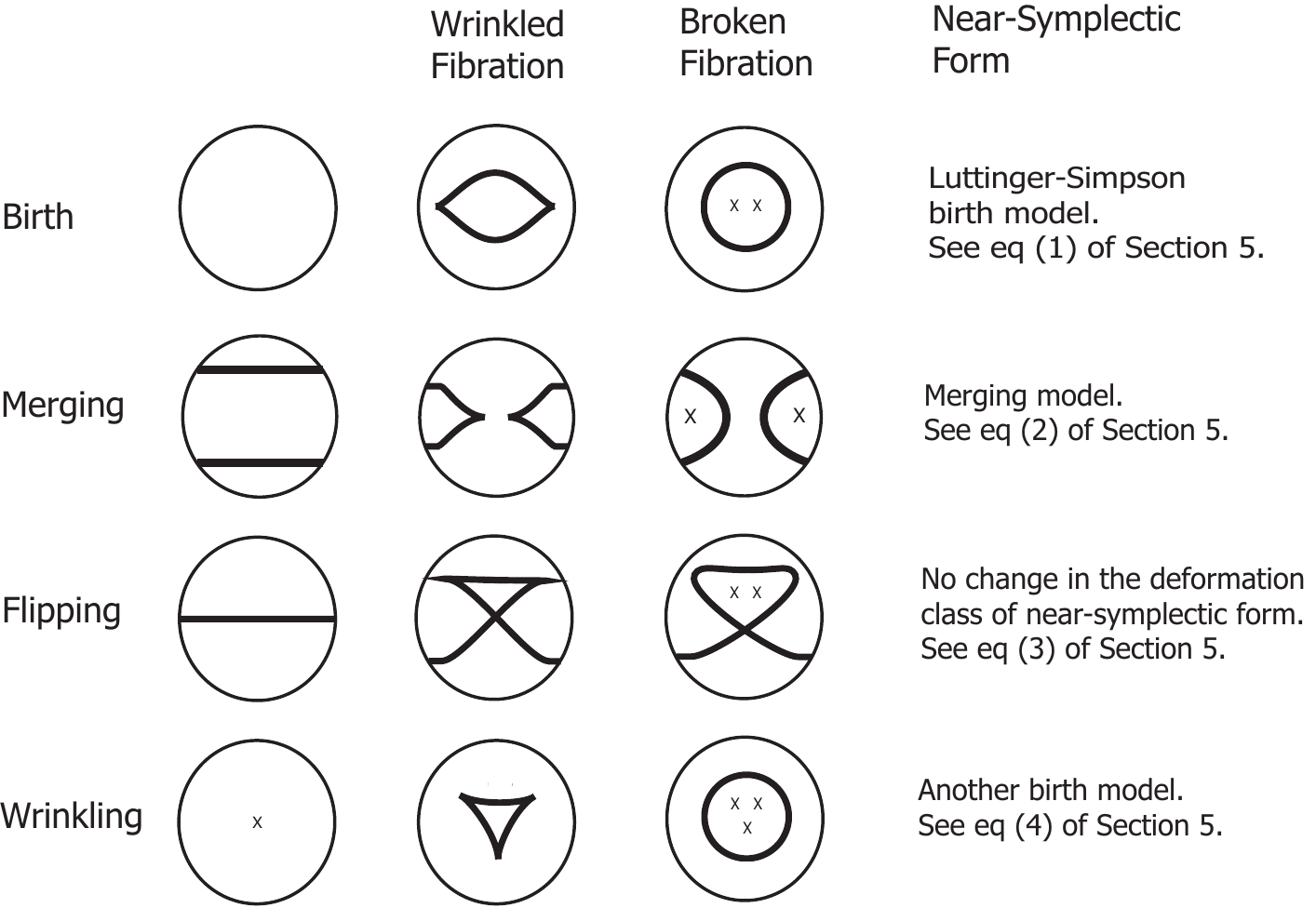}
\end{center}
\caption{Table of Moves}	 
\end{figure} 

The next important task that we would like to address in the future is to prove that the Lagrangian matching invariant that was described in the introduction is invariant under the set of moves described in this paper. Of equal importance is the problem of determining the set of equivalence classes of deformations of broken fibrations on a given $4$--manifold. The author believes that homotopic broken fibrations should be deformation equivalent. That is, we would like to prove that some sort of \emph{h}--principle holds for wrinkled fibrations. The main difficulty here is that wrinkled fibrations are constrained to have indefinite Hessian along the critical points. 
 
\nocite{yasha}

\newpage
\appendix
\section{Classification of $(1,1)$--stable unfoldings}
\label{appendix}
Throughout, we denote by $\mathcal{E}(n)$ the set of germs at $0\in \f{R}^n$ of smooth mappings from $\f{R}^n$ to $\f{R}$, and $\mathfrak{m}(n)$ is the unique maximal ideal of $\mathcal{E}(n)$ consisting of those germs $f$ such that $f(0)=0$.
\begin{definition}
Let $\eta\in \mathfrak{m}(n)$. An $r$--dimensional unfolding of $\eta$ is a germ $f\in \mathcal{E}(n+r)$ such that $f|\f{R}^n = \eta$.
\end{definition}

For the next definition, recall the definition of $(1,1)$--equivalence of unfoldings was given in Definition \ref{defrsequiv}.

\begin{definition}
	
Let $f\in \mathfrak{m}(n+d+2)$ and let $g\in \mathfrak{m}(n+2)$. We say $f$ $(1,1)$--reduces to $g$ if there is a non-degenerate quadratic form $Q$ on $\f{R}^d$ such that $f$ is $(1,1)$--equivalent to the germ $g' \in \mathfrak{m}(n+d+2)$ given by $g'(s,t,x,y) = g(s,t,x)+Q(y)$ for $s\in \f{R}, t\in \f{R},x\in \f{R}^n, y\in \f{R}^d$.
\end{definition}

We will give a classification $(1,1)$--stable unfoldings up to $(1,1)$--equivalence based on the algorithm described in \cite{wasserman}. We shall make direct use of the lemmas and theorems in \cite{wasserman} without restating them here. 

\begin{theorem}
\label{theoremapp}	
Let $f \in \mathfrak{m}(n+2)$ be a $(1,1)$--stable unfolding of $\eta \in \mathfrak{m}(n)^2$. Then either $f$ has Morse singularity at $0$, or $f$ $(1,1)$--reduces to a unique one of the following unfoldings $h_i$ of germs $v_i$ :
\begin{eqnarray*} 
  v_i \ \ \ \ \ \ \  &  &  \ \ \ \ \  h_i \ \ \ \ \ \ \ \ \ \ \\ 
  v_0(x)=x^3 &   & h_0(s,t,x)=x^3+tx \\
  v_1(x)=x^3 &   & h_1(s,t,x)=x^3+t^2x+sx \\ 
  v_2(x)=x^3 &   & h_2(s,t,x)=x^3-t^2x+sx \\
  v_3(x)=x^4 &   & h_3(s,t,x)=x^4+sx^2+tx 
\end{eqnarray*}
\end{theorem}

We follow the same method as Wasserman's classification of $(3,1)$--stable unfoldings in \cite{wasserman}. In particular, the following special case of Theorem $4.11$ from \cite{wasserman} plays a crucial role in the classification. We say that $f\in \mathfrak{m}(n+2)$ is $2$--stable if \[ \mathcal{E}(u,x)=\langle \frac{\del f}{\del x} \rangle \mathcal{E}(u,x) +\langle \frac{\del f}{\del u} \rangle \mathcal{E}(u)+F^*\mathcal{E}(\f{R}^3) \] 
where $F(u,x)=(u,f(u,x))$ for $u\in \f{R}^2, x\in \f{R}^n$. This notion is very similar to $(1,1)$--stability which was given in Definition \ref{defstable}. The difference is that here we do not distinguish the variables $s$ and $t$. In particular, $(1,1)$--stability implies $2$--stability.

\begin{theorem}
\label{theoremwass}	
Let $g \in \mathfrak{m}(n+2)$ be a $(1,1)$--stable unfolding of $\eta\in \mathfrak{m}(n)$, and suppose that $f\in \mathfrak{m}(n+2)$ is a $2$--stable unfolding of $\eta$. Then there exists a polynomial germ $p \in \mathfrak{m}(\f{R})$ of degree at most $2$ such that $g$ is $(1,1)$--equivalent to either $f(s+p(t),t, x)$ or $f(t,s+p(t),x)$.	
	
\end{theorem}

\emph{Proof of Theorem \ref{theoremapp}:} If $f$ is $(1,1)$--stable, then $f$ is $2$--stable and hence $f$ has a simple singularity or $f$ reduces (in the sense of Definition $2.24$ of \cite{wasserman}) to a unique one of the unfoldings $g_i$ in Thom's list of seven elementary catastrophes (see Theorem $2.20$ of \cite{wasserman}). If the latter case occurs, then $\eta$ reduces to a unique one of the germs $\mu_i$ in Thom's list. By Lemma $4.18$ of \cite{wasserman}, if $\eta$ reduces to $\mu_i$ then $f$ $(1,1)$--reduces to a two-dimensional unfolding $h$ of $\mu_i$ which by Lemma $4.17$ of \cite{wasserman} must be $(1,1)$--stable. Moreover, Lemma $4.19$ together with Lemma $4.20$ in \cite{wasserman} implies that the set of $(1,1)$--stable unfoldings of $\mu_i$ to which $f$ $(1,1)$--reduces is exactly the $(1,1)$--equivalence class of $h$. Hence to complete the proof we need only to show that for each germ $\mu_i$ in Thom's list, the list of Theorem \ref{theoremapp} gives exactly the classification of $(1,1)$--stable unfoldings of $\mu_i$ up to $(1,1)$--equivalence.

First, consider the case of $\mu_1(x)=x^3.$ By Theorem \ref{theoremwass}, a $(1,1)$--stable unfolding $f$ up to $(1,1)$--equivalence of $\mu_1(x)=x^3$ is either $x^3+tx$ or of the form $x^3+(s+at^2+bt)x$ where $a, b \in \f{R}$. The former case is $h_0$, so we concentrate on the latter case.

Corollary $4.13$ of \cite{wasserman} gives the $(1,1)$--stable condition for $f$ as :
\[ \mathcal{E}(t,x) = \langle \frac{\del f_0}{\del x} \rangle {\mathcal{E}(t,x)} + \langle \frac{\del f_0 } {\del t} \rangle {\mathcal{E}(t)} +\f{R}\langle \frac{\del f}{\del s}|_{\{s=0\}} \rangle + \langle 1, f_0 \rangle {\mathcal{E}(t)} +  \mathfrak{m}(t)^2 \mathcal{E}(t,x)+\mathfrak{m}(t,x)^4  \]
where $f_0 = f|_{\{s=0\}}$. Thus, $f$ is $(1,1)$--stable if and only if 
\begin{eqnarray*} 
\mathcal{E}(t,x) & = & \langle 3x^2+(at^2+bt)\rangle {\mathcal{E}(t,x)} + \langle 2atx+bx \rangle {\mathcal{E}(t)} +\f{R}\langle x \rangle + \langle 1, x^3+(at^2+bt)x \rangle {\mathcal{E}(t)} \\
& + & \mathfrak{m}(t)^2 \mathcal{E}(t,x)+\mathfrak{m}(t,x)^4.
\end{eqnarray*}
An easy calculation then reveals that $f$ is $(1,1)$--stable if and only if $a$ or $b$ is nonzero. Suppose $b\neq 0$, then we change coordinates by setting $t'=s+at^2+bt, s'=s$ and get $f$ is $(1,1)$--equivalent to $h_0$. On the other hand, if $b=0$, then by scaling $t$, we obtain that $f$ is $(1,1)$--equivalent to either $h_1$ or $h_2$. Furthermore, it is clear that none of $h_0$, $h_1$ and $h_2$ are $(1,1)$--equivalent.

Now consider the case of $\mu_2(x)=x^4$. By Theorem \ref{theoremwass}, a $(1,1)$--stable unfolding $f$ up to $(1,1)$--equivalence of $\mu_2(x)=x^4$ is either $x^4+(s+at^2+bt)x^2+tx$ or $x^4+tx^2+(s+at^2+bt)x$, where $a,b \in \f{R}$. In order to determine for which values of $a$ and $b$ these maps are $(1,1)$--stable, we again apply the criteria given by Corollary $4.13$ of \cite{wasserman}. Suppose first $f$ is given by $x^4+(s+at^2+bt)x^2+tx$. Then $f$ is $(1,1)$--stable if and only if 
\begin{eqnarray*} 
\mathcal{E}(t,x) & = & \langle 4x^3+2at^2x+2btx+t \rangle {\mathcal{E}(t,x)} + \langle 2atx^2+bx^2+x   \rangle {\mathcal{E}(t)} +\f{R}\langle x^2 \rangle \\
& + & \langle 1, x^4+(at^2+bt)x^2+tx \rangle {\mathcal{E}(t)} +\mathfrak{m}(t)^2\mathcal{E}(t,x)+\mathfrak{m}(t,x)^4.
\end{eqnarray*}
It turns out that in this case $f$ is $(1,1)$--stable for all values of $a$ and $b$. By Lemma $4.9$ of \cite{wasserman} stably homotopic $(1,1)$--stable germs are $(1,1)$--equivalent. Therefore we can set $a=b=0$ and conclude that $f$ is $(1,1)$--equivalent to $h_3=x^4+sx^2+tx$.
 
Finally, suppose that $f$ is given by $x^4+tx^2+(s+at^2+bt)x$. Then Corollary $4.13$ of \cite{wasserman} yields that $f$ is $(1,1)$--stable if and only if
\begin{eqnarray*} 
\mathcal{E}(t,x)&=&\langle 4x^3+2tx+at^2+bt \rangle {\mathcal{E}(t,x)} + \langle x^2+2atx+bx   \rangle {\mathcal{E}(t)} +\f{R}\langle x \rangle \\
& + & \langle 1, x^4+tx^2+(at^2+bt)x \rangle {\mathcal{E}(t)} +\mathfrak{m}(t)^2 \mathcal{E}(t,x)+\mathfrak{m}(t,x)^4.
\end{eqnarray*}
It is then an easy calculation to conclude that $f$ is $(1,1)$--stable if and only if $b$ is nonzero. Next, we again apply Lemma $4.9$ of \cite{wasserman} to set $a=0$, and conclude that $f$ is $(1,1)$--equivalent to $x^4+tx^2+(s+bt)x$. Now, we change coordinates by setting $t'=s+bt, s'=-s/b$. Then $f$ becomes $x^4+(s'+t'/b)x^2+t'x$. So we are back to the previous case, hence we conclude that $f$ is $(1,1)$--equivalent to $h_3$, as desired.

\end{document}